\documentclass[12pt]{amsart}
\usepackage{amsmath}
\usepackage{amssymb}
\usepackage{amsfonts}

\newtheorem{theorem}{Theorem}[section]
\theoremstyle{plain}

\newtheorem{definition}[theorem]{Definition}

\newtheorem{lemma}[theorem]{Lemma}

\newtheorem{proposition}[theorem]{Proposition}
\newtheorem{remark}[theorem]{Remark}

\numberwithin{equation}{section}

\begin{document}
\title[Crossed-products for conformal automorphisms of $\mathbb{D}$]{Crossed-product C*-algebras for conformal automorphisms of the disk}
\author{Man-Duen Choi}
\address{Department of Mathematics, University of Toronto}
\email{choi@math.toronto.edu.}
\author{Fr\'{e}d\'{e}ric Latr\'{e}moli\`{e}re}
\address{Department of Mathematics, University of Toronto}
\email{frederic@math.berkeley.edu}
\urladdr{http://www.math.toronto.edu/frederic/}

\begin{abstract}
We study the C*-algebra crossed-product of the closed unit disk by the action
of one of its conformal automorphisms. After classifying the conformal
automorphisms up to topological conjugacy, we investigate, for each class, the
irreducible representations of the full C*-crossed-products, and derive their
spectrum and a complete desciption of the algebras.

\end{abstract}
\maketitle

\section{Introduction}

We present in this paper a detailed description of the C*-crossed product
algebras associated to conformal discrete dynamical systems of the closed unit
disk $\mathbb{D}$.

Given a compact Hausdorff space $X$ and an homeomorphism $\varphi$ on $X$, the
C*-crossed product of $C(X)$ by $\mathbb{Z}$ generated by the *-automorphism
$f\in C(X)\mapsto f\circ\varphi$ is the universal C*-algebra
\cite{Zeller-Meier68} \cite{Pedersen79} \cite{Tomiyama87} generated by $C(X)$
and a unitary $U_{\varphi}$ such that for all $f\in C(X)$ we have $U_{\varphi
}^{\ast}fU_{\varphi}=f\circ\varphi$. We refer to $U_{\varphi}$ as the
canonical unitary of $C(X)\rtimes_{\varphi}\mathbb{Z}$. The universal property
means that for any C*-algebra $\mathfrak{A}$ containing a unitary $V$ and a
C*-algebra $\mathcal{A}$ such that there exists a *-homomorphism
$\rho:C(X)\longrightarrow\mathcal{A}$ which satisfies $\rho(U_{\varphi}^{\ast
}fU_{\varphi})=V^{\ast}\rho(f)V$ for all $f\in C(X)$, there exists a
(necessarily unique)\ *-homomorphism $\psi:C(X)\rtimes_{\varphi}%
\mathbb{Z}\longrightarrow\mathfrak{A}$ such that $\psi$ restricts to $\rho$ on
$C(X)$ and $\psi(U_{\varphi})=V$. A construction of this C*-algebra can be
found in \cite{Zeller-Meier68} \cite{Pedersen79}. Given two compact separable
spaces $X$ and $Y$ and two homeomorphisms $\varphi$ and $\psi$ of $X$ and $Y$
respectively, if $\varphi$ and $\psi$ are topologically conjugated in the
sense that there exists an homeomorphism $\gamma$ from $X$ onto $Y$ such that
$\gamma\circ\varphi\circ\gamma^{-1}=\psi$, then by the universal property we
have $C(X)\rtimes_{\varphi}\mathbb{Z=}C(Y)\rtimes_{\psi}\mathbb{Z}$. Thus,
C*-crossed-products are invariants of topological conjugacy, though they are
known not to be complete invariants.

We study in this paper the case when $X=\mathbb{D}=\left\{  z\in
\mathbb{C}:\left\vert z\right\vert \leq1\right\}  $ and $\varphi$ is a
conformal automorphism of $\mathbb{D}$. The C*-algebra $C(\mathbb{D})$ can be
represented as the C*-algebra $C^{\ast}(A_{\varphi})$ for any normal element
$A_{\varphi}$ on some Hilbert space with $\sigma(A_{\varphi})=\mathbb{D}$
where $\sigma(.)$ denotes the spectrum of an operator. We denote the canonical
unitary of $C(\mathbb{D})\rtimes_{\varphi}\mathbb{Z}$ by $U_{\varphi}$. Thus,
$C(\mathbb{D})\rtimes_{\varphi}\mathbb{Z}$ is the universal C*-algebra
generated by a unitary $U_{\varphi}$ and a normal element $A_{\varphi}$ such
that $\sigma(A_{\varphi})=\mathbb{D}$ and $U_{\varphi}^{\ast}A_{\varphi
}U_{\varphi}=\varphi(A_{\varphi})$. Any pair $(A,U)$ where $A$ is a normal
operator with $\sigma(A)\subseteq\mathbb{D}$ and $U$ a unitary operator on
some Hilbert space with $U^{\ast}AU=\varphi(A)$ will be called a covariant
pair for $\left(  C(\mathbb{D}),\varphi,\mathbb{Z}\right)  $. Since
$\sigma(A)$ is a compact subset of $\sigma(A_{\varphi})=\mathbb{D}$, we define
a *-epimorphism $\rho$ from $C\left(  \mathbb{D}\right)  $ onto $C\left(
\sigma(A)\right)  $ by setting $\rho(f(A_{\varphi}))=f(A)$ for any $f\in
C(\mathbb{D})$. We observe that $U^{\ast}f(A)U=f\circ\varphi(A)=f\left(
U_{\varphi}^{\ast}A_{\varphi}U_{\varphi}\right)  $, so by universality there
exists a unique *-epimorphism $\psi:C(\mathbb{D})\rtimes_{\varphi}%
\mathbb{Z}\twoheadrightarrow C^{\ast}(A,U)$ such that $\psi(A_{\varphi})=A$
and $\psi(U_{\varphi})=U$.

While C*-crossed products of minimal dynamical systems have received a lot of
attention, more general and natural situations have not been studied as much.
For crossed-product for action on the circle, one can consult \cite{Rango02}.
We offer in this paper to investigate a few concrete examples of non-proper,
non-free, non-minimal dynamical systems. Some of our examples are of type I,
and some general properties of such crossed-products were derived in
\cite{Tomiyama02}.

\bigskip We start this paper with a complete classification of the conformal
automorphisms of $\mathbb{D}$ up to topological conjugacy. We split the group
$\mathbb{M}_{\mathbb{D}}$ of conformal automorphisms of $\mathbb{D}$ into the
set of hyperbolic automorphisms, which are all topologically conjugated with
one another, the set of parabolic automorphisms which are all conformally
conjugated with each other, the set of elliptic automorphisms, each of which
is conjugated to a nontrivial rotation around the origin, and the identity
map, which is a special case whose crossed-product is trivially $C(\mathbb{D}%
)\rtimes_{\operatorname*{Id}}\mathbb{Z}=C(\mathbb{D}\times\mathbb{T})$.

We then study the full C*-crossed-product algebra for each conformal
automorphism. In each case, we provide a full description of its Gelfan'd
spectrum, with its Zarisky (or hull-kernel) topology. We then describe all the
irreducible representations, explicitly for the hyperbolic and parabolic
cases, or based upon the representation theory of the rotation algebras for
the elliptic case. Last, we derive an explicit description of each
crossed-product in term of more elementary algebras, exposing the topological
content of the C*-crossed-product in each case. We prove that:

\begin{itemize}
\item For hyperbolic automorphisms, the crossed-product is isomorphic to
\[
\left\{
\begin{array}
[c]{c}%
f\in C\left(  \mathbb{R}\times\lbrack0,1],\left[
\begin{array}
[c]{cc}%
\mathcal{T} & \mathcal{K}\\
\mathcal{K} & \mathcal{T}%
\end{array}
\right]  \right)  \text{ such that:}\\
\forall x\in\mathbb{R}\ \forall y\in\lbrack0,1]\ \ \left\{
\begin{array}
[c]{c}%
f(x,y)-f(0,0)\in\mathcal{K}\\
U^{\ast}f(x,y)U=f(x+1,y)
\end{array}
\right.
\end{array}
\right\}  ,
\]
where $\mathcal{T}$ is the Toeplitz algebra and $\mathcal{K}$ is the compact
operator algebra,

\item For parabolic automorphisms, the crossed-product is isomorphic to
$C(U^{\prime})+C_{1}(\mathbb{D},\mathcal{K})$ where $U^{\prime}:t\in
\mathbb{D}\mapsto U$ with $U$ the bilateral shift on $l^{2}\left(
\mathbb{Z}\right)  $ and $\mathcal{K}$ the compact operator algebras acting on
$l^{2}\left(  \mathbb{Z}\right)  $,

\item For an elliptic automorphism conjugated to the rotation of angle
$2\pi\theta$, the crossed-product is isomorphic to
\[
\left\{  f\in C\left(  [0,1],\mathcal{A}_{\theta}\right)  :f(0)\in C^{\ast
}(U)\right\}  \text{,}%
\]
where $\mathcal{A}_{\theta}=C^{\ast}(U,V)$ is the universal C*-algebra
generated by two unitaries $U,V$ such that $VU=e^{2i\pi\theta}UV$.
\end{itemize}

In this paper, unless otherwise stated, the norm of a Banach space $E$ is
denoted by $\left\Vert .\right\Vert _{E}$.

\section{Classification of the conformal automorphisms of the closed unit
disk}

The following well-known corollary of the Schwarz lemma can be found in many
textbooks in complex analysis (e.g. \cite[Theorem 2.1 p.213]{Lang99}):

\begin{proposition}
\label{autodisk}Let $\varphi$ be a conformal automorphism of the closed unit
disk $\mathbb{D=}\left\{  z\in\mathbb{C}:\left\vert z\right\vert
\leq1\right\}  $. Then there exists $z_{0}\in\mathbb{C}$ and $\theta\in
\lbrack0,1)$ such that $\left\vert z_{0}\right\vert <1$ and:%
\[
\varphi(z)=\exp\left(  2i\pi\theta\right)  \frac{z-z_{0}}{1-\overline{z_{0}}z}%
\]
for all $z\in\mathbb{D}$.

Moreover, $\varphi(0)=0$ if and only if there exists $\theta\in\lbrack0,1)$
such that for all $z\in\mathbb{D}$ we have $\varphi(z)=\exp\left(  2i\pi
\theta\right)  z$.
\end{proposition}

We recall the following well-known definitions:

\begin{definition}
Let $\psi$ be an homeomorphism of $\mathbb{D}$. A point $\omega\in\mathbb{D}$
is attractive for the dynamical system $\left(  C(\mathbb{D}),\psi
,\mathbb{Z}\right)  $ when there exists an open neighborhood $\Omega$ of
$\omega$ in $\mathbb{D}$ such that for all $y\in\Omega$ we have:
$\lim_{n\rightarrow\infty}\psi^{n}(y)=\omega$. A point $\alpha$ is repulsive
when there exists an open neighborhood $\Omega$ of $\alpha$ in $\mathbb{D}$
such that $\lim_{n\rightarrow\infty}\psi^{-n}(y)=\alpha$ for all $y\in\Omega$.
\end{definition}

We now classify all the conformal automorphisms of $\mathbb{D}$ up to
topological conjugacy, since C*-crossed-products only depend on topological
conjugacy classes. As a first step, we identify the conformal conjugacy
classes of conformal automorphisms of $\mathbb{D}$.

\begin{theorem}
\label{autodisk2}Let $\varphi\in\mathbb{M}_{\mathbb{D}}$ where $\mathbb{M}%
_{\mathbb{D}}$ is the group of conformal automorphisms of $\mathbb{D}$: thus
$\varphi:z\in\mathbb{D}\mapsto\exp\left(  2i\pi\theta\right)  \frac{z-z_{0}%
}{1-\overline{z_{0}}z}$ for some $\theta\in\lbrack0,1)$ and $z_{0}%
\in\mathbb{C}$ such that $\left\vert z_{0}\right\vert <1$. Then:

\begin{enumerate}
\item The following are equivalent:

\begin{itemize}
\item $\left\vert z_{0}\right\vert >\left\vert \sin\left(  \pi\theta\right)
\right\vert $,

\item $\varphi$ has exactly two distinct fixed points in $\mathbb{D}$,

\item $\varphi$ has exactly two fixed points on $\mathbb{T}$,

\item $\varphi$ is conjugated in $\mathbb{M}_{\mathbb{D}}$ with $z\in
\mathbb{D}\mapsto\frac{z+a}{1+az}$ for some unique $a\in(0,1)$.
\end{itemize}

If these assertions hold, then $\varphi$ is called \emph{hyperbolic} and one
fixed point of $\varphi$ is attractive, while the other one is repulsive. In
particular, $z\in\mathbb{D}\mapsto\frac{z+a}{1+az}$ for $a\in(0,1)$ admits $1$
as attractive fixed point and $-1$ as repulsive fixed point.

\item The following are equivalent:

\begin{itemize}
\item $\left\vert \sin\left(  \pi\theta\right)  \right\vert >\left\vert
z_{0}\right\vert $,

\item $\varphi$ has exactly one fixed point $\alpha$ in $\mathbb{D}$ such that
$\left\vert \alpha\right\vert <1$,

\item $\varphi$ is conjugated in $\mathbb{M}_{\mathbb{D}}$ to a rotation
$z\in\mathbb{D}\mapsto\mu z$ for some unique $\mu\in\mathbb{T}\backslash\{1\}$.
\end{itemize}

If these assertions hold, then $\varphi$ is called \emph{elliptic}.

\item The following are equivalent:

\begin{itemize}
\item $\left\vert \sin\left(  \pi\theta\right)  \right\vert =\left\vert
z_{0}\right\vert \not =0$,

\item $\varphi$ has exactly one fixed point $\alpha$ in $\mathbb{D}$ such that
$\left\vert \alpha\right\vert =1$,

\item if $\varphi(-1)$ is on the upper half of $\mathbb{T}$, then $\varphi$ is
conjugated in $\mathbb{M}_{\mathbb{D}}$ with $\phi:z\mapsto\frac{3i-1}%
{i-3}\frac{z-\frac{2i-1}{i-3}}{1-\frac{-2i-1}{-i-3}z}$ , where $\phi$ is the
unique conformal automorphism of $\mathbb{D}$ whose only fixed point is $1$
and $\phi(-1)=i$.

\item If $\varphi(-1)$ is on the lower half of $\mathbb{T}$, then $\varphi$ is
conjugated in $\mathbb{M}_{\mathbb{D}}$ with $\phi:z\mapsto\frac{3i+1}%
{i+3}\frac{z-\frac{2i+1}{i+3}}{1+\frac{2i-1}{-i+3}z}$ , where $\phi$ is the
unique conformal automorphism of $\mathbb{D}$ whose only fixed point is $1$
and $\phi(-1)=-i$.
\end{itemize}

If these assertions hold, $\varphi$ is called \emph{parabolic} and its fixed
point is both attractive and repulsive.

\item Last, $\varphi$ is the identity map if and only if it has more than two
fixed points, if and only if $\theta=0$ and $z_{0}=0$, if and only if it has
two fixed points $\alpha,\beta$ with either $\left\vert \alpha\right\vert
=1,\left\vert \beta\right\vert <1$ or $\left\vert \alpha\right\vert
<1,\left\vert \beta\right\vert <1$.
\end{enumerate}
\end{theorem}

\begin{proof}
The fixed points of $\varphi$ are the solutions of $\overline{z_{0}}%
z^{2}+(e^{2i\pi\theta}-1)z-e^{2i\pi\theta}z_{0}=0$. Thus, by a straightforward
computation, the location of the fixed points of $\varphi\not =%
\operatorname*{Id}$ are determined by the sign of $\left\vert z_{0}\right\vert
-\left\vert \sin\left(  \pi\theta\right)  \right\vert $ as stated in the
theorem (\textit{see appendix at the end for the calculation}).

We now turn to the matter of proving that $\varphi$ is, in each case,
conjugated to a specific conformal automorphism. First, let us assume that
$\varphi$ has two distinct fixed points $\alpha,\beta\in\mathbb{T}$. Assume
first $\alpha,\beta\not \in \left\{  -1,1\right\}  $. Let $\omega$ be a square
root of $\overline{\alpha\beta}$ and let $\gamma=\frac{i-\omega\alpha
}{1-i\omega\alpha}$ (observe that $\gamma\in]-1,1[$). We set $\psi
(z)=i\omega\frac{z+\gamma\overline{\omega}}{1+\gamma\omega z}$ for all
$z\in\mathbb{D}$. Then $\psi\circ\varphi\circ\psi^{-1}$ fixes $-1$ and $1$.

Now, any conformal automorphism $\varphi$ has $-1$ and $1$ as fixed points if
and only if $\lambda-1=0$ and $z_{0}/\overline{z_{0}}=1$, namely $z_{0}%
\in\mathbb{R}$ and $e^{2i\pi\theta}=1$. Up to conjugating $\varphi$ by
$\psi^{\prime}:z\mapsto-z$ we can assume $z_{0}\leq0$, so all hyperbolic
automorphisms of $\mathbb{D}$ are conjugated in $\mathbb{M}_{\mathbb{D}}$ to
$z\mapsto\frac{z+z_{0}}{1+z_{0}z}$ for some $z_{0}\in(0,1)$, which a simple
calculation shows to be uniquely determined by $\varphi$ (\textit{see appendix
at the end for the calculation}).

\qquad If instead $\left\{  \alpha,\beta\right\}  \not \subset \mathbb{T}$
then $\varphi$ has only one fixed point $c\in\{\alpha,\beta\}$ in $\mathbb{D}%
$. Then set $\psi:z\mapsto\frac{z-c}{1-\overline{c}z}\in\mathbb{M}%
_{\mathbb{D}}$. The map $R=\psi\circ\varphi\circ\psi^{-1}$ is a conformal
automorphism of $\mathbb{D}$ with $R(0)=0$. Hence $R$ is a rotation by Theorem
(\ref{autodisk}). Now, if $\varphi\in\mathbb{M}_{\mathbb{D}}$ conjugates
$z\mapsto\mu z$ with $z\mapsto\nu z$ for $\mu,\nu\in\mathbb{T}$ then
$\varphi(0)=0$ so $\varphi:z\mapsto\eta z$ for $\eta\in\mathbb{T}$ and thus
$\mu=\nu$.

\qquad Last, observe that when $\varphi$ is parabolic with fixed point
$\alpha$ then $r_{\alpha^{-1}}\circ\varphi\circ r_{\alpha}(1)=1$, where
$r_{\alpha}:z\mapsto\alpha z$. Hence up to conformal conjugacy, we may as well
assume that $\varphi$ admits $1$ as its unique fixed point. Let $z_{0}$ with
$\left\vert z_{0}\right\vert <1$ and $\lambda\in\mathbb{T}$ such that
$\varphi:z\mapsto\lambda\frac{z-z_{0}}{1-\overline{z_{0}}z}$. Now
$\varphi(1)=1$ if and only if $1=\frac{1-\lambda}{2\overline{z_{0}}}$, so if
and only if $\varphi:z\mapsto\lambda\frac{z-\frac{1-\overline{\lambda}}{2}%
}{1-\frac{1-\lambda}{2}z}$. Now, if $\varphi:z\mapsto\lambda\frac
{z-\frac{1-\overline{\lambda}}{2}}{1-\frac{1-\lambda}{2}z}$ for some
$\lambda\in\mathbb{T}$, then let $\psi$ be the (unique) conformal automorphism
of $\mathbb{D}$ such that $\psi(1)=1$, $\psi(-1)=-1$ and $\psi(\varphi\left(
-i\right)  )=i$ if $\varphi(-i)$ is in the upper half plane, or $\psi
(\varphi(-1))=-i$ otherwise. Then $\phi=\psi\circ\varphi\circ\psi^{-1}$ is
parabolic since it has a unique fixed point in $\mathbb{T}$ and we have
$\phi(1)=1$ and $\phi(-1)=i$ if $\varphi(-1)$ is in the upper half plane, and
$\phi(-1)=-i$ otherwise. We deduce that $\phi:z\mapsto\frac{3i-1}{i-3}%
\frac{z-\frac{2i-1}{i-3}}{1-\frac{-2i-1}{-i-3}z}$ (if $\varphi(-1)$ is in the
upper half-plane) or $\phi:z\mapsto\frac{3i+1}{i+3}\frac{z-\frac{2i+1}{i+3}%
}{1+\frac{2i-1}{-i+3}z}$. This concludes our proof.
\end{proof}

Now, we have reduced the topological conjugacy classification problem to
classifying the special automorphisms in Proposition (\ref{autodisk2}):

\begin{theorem}
\label{topoconj}Let $\varphi,\psi$ be two conformal automorphisms of
$\mathbb{D}$. Then:

\begin{itemize}
\item If $\varphi$ and $\psi$ are topologically conjugated and if $\varphi$ is
hyperbolic (resp. parabolic, elliptic) then $\psi$ is hyperbolic (resp.
parabolic, elliptic),

\item If $\varphi$ is elliptic then $\varphi$ is topologically conjugated to a
rotation $r_{\theta}:z\in\mathbb{D}\mapsto e^{2i\pi\theta}z$ for some
$\theta\in(0,1)$. Moreover, $r_{\theta}$ and $r_{\theta^{\prime}}$ are
topologically conjugated if and only if $\theta\in\left\{  \theta^{\prime
},-\theta^{\prime}\right\}  $.

\item All parabolic automorphisms are topologically conjugated on $\mathbb{D}
$,

\item All hyperbolic automorphisms are topologically conjugated on
$\mathbb{D}$. Moreover, if $\varphi:z\mapsto\frac{z+a}{1+az}$ for some
$a\in(0,1)$ and if we set: $$\mathcal{D}_{\varphi}=\left\{  z\in\mathbb{D}%
:\operatorname{Re}(z)\geq0,\left\vert 1-z\right\vert \geq1-a\right\}$$ then 
for any $x \in \mathbb{D}\backslash\left\{  -1,1\right\}  $ there exists a
unique $y\in\mathcal{D}_{\varphi}$ such that $\mathcal{O}_{\varphi
}(x)=\mathcal{O}_{\varphi}(y)$, where $\mathcal{O}_{\varphi}(w)=\left\{
\varphi^{n}(w):n\in\mathbb{Z}\right\}  $ is the orbit for $\varphi$ of any
$w\in\mathbb{D}$.
\end{itemize}
\end{theorem}

\begin{proof}
It is well known that two rotations $r_{\lambda}:z\mapsto\lambda z$ and
$r_{\mu}:z\mapsto\mu z$ for $\lambda,\mu\in\mathbb{T}$ are conjugated if and
only if $\lambda\in\left\{  \mu,\overline{\mu}\right\}  $: if $\lambda=\mu$ it
is trivial, if $\lambda=\overline{\mu}$ then $\overline{\left(  r_{\lambda
}(\overline{z})\right)  }=r_{\mu}$, so the condition is sufficient, and it is
necessary by \cite[Theorem 11.2.7, p. 397]{Katok95}. Since the number of fixed
points is a topological invariant, and any homeomorphism of $\mathbb{D}$ maps
$\mathbb{T}$ to $\mathbb{T}$, the two first statements of our theorem are proven.

By Theorem (\ref{autodisk2}), if $\varphi$ is parabolic then it is
conformally, hence topologically conjugated to either \[\phi_{+}:z\mapsto
\frac{3i-1}{i-3}\frac{z-\frac{2i-1}{i-3}}{1-\frac{-2i-1}{-i-3}z} \;\; \textrm{or} \;\; \phi
_{-}:z\mapsto\frac{3i+1}{i+3}\frac{z-\frac{2i+1}{i+3}}{1+\frac{2i-1}{-i+3}z}\textrm{.}\]
Yet, if $\nu:z\mapsto\overline{z}$ then $\nu^{-1}\circ\phi_{-}\circ\nu
=\phi_{+}$, hence the third assertion of our theorem.

By Proposition (\ref{autodisk2}), any hyperbolic automorphism $f$ is
conformally, hence topologically conjugated in $\mathbb{D}$ to an hyperbolic
automorphism $\varphi$ with fixed points $-1$ and $1$ where $1$ is the
attractive fixed point. It suffices therefore to prove that any two hyperbolic
automorphisms $\varphi$ and $\psi$ of $\mathbb{D}$ with fixed points $-1$ and
$1$, where $1$ is attractive, are topologically conjugated. Therefore, let
$\varphi:z\mapsto\frac{z+a}{1+az}$ and $\psi:z\mapsto\frac{z+b}{1+bz}$ for
$a,b\in(0,1)$.

We define the region $\overline{\mathcal{D}_{\varphi}}$ as the closed
subset of $\mathbb{D}$ contained between the diameter $L_{\varphi}=\left\{
it:t\in\lbrack-1,1]\right\}  $ and its image $\varphi(L)$, which is the circle
of center $\varphi\left(  \infty\right)  =\frac{1}{a}$ and radius $\frac{1}
{a}-a$. We denote by $\mathcal{D}_{\varphi}$ the set $\overline{\mathcal{D}
_{\varphi}}\backslash\varphi(L)$ and we check easily that the orbit for
$\varphi$ of any point of $\mathbb{D}\backslash\left\{  -1,1\right\}  $
intersects $\mathcal{D}_{\varphi}$ at exactly one point.

The region $\overline{\mathcal{D}_{\varphi}}$ is homeomorphic to the rectangle
$[0,1]^{2}$. We choose an homeomorphism $\mu_{\varphi}:\overline
{\mathcal{D}_{\varphi}}\rightarrow\lbrack0,1]^{2}$ such that $\mu
_{\varphi}(it)=(t,0)$ and $\mu_{\varphi}\left(  \varphi(it\right)  )=(t,1)$
for all $t\in\lbrack-1,1]$.

Of course, we can as well construct an homeomorphism $\mu_{\psi}
:\overline{\mathcal{D}_{\psi}}\longrightarrow\lbrack0,1]^{2}$ such that
$\mu_{\psi}(it)=(t,0)$ and $\mu_{\psi}(\psi(it))=\left(  t,1\right)  $. Hence
$\mu=\mu_{\psi}^{-1}\mu_{\varphi}$ is an homeomophism from $\overline
{\mathcal{D}_{\varphi}}$ onto $\overline{\mathcal{D}_{\psi}}$ such that
$\mu\circ\varphi=\psi\circ\mu$ on the set $\left\{  it:t\in\lbrack
-1,1]\right\}  $. By a trivial induction, we can extends $\mu$ to the set
$\mathbb{D}\backslash\left\{  -1,1\right\}  $ and then set $\mu(-1)=-1$ and
$\mu(1)=1$, and this map is easily checked to be an homeomorphism of
$\mathbb{D}$ such that $\mu\circ\varphi=\psi\circ\mu$ on $\mathbb{D}$ (\textit{see appendix for a more explicit construction of an example of }$\mu$).
\end{proof}

This achieves the classification, up to topological classification, of all the
conformal automorphisms of $\mathbb{D}$. In particular, the description of the
dynamical system $\left(  C(\mathbb{D}),\varphi,\mathbb{Z}\right)  $ is now
easy. If $\varphi$ is hyperbolic, then (up to conformal conjugacy) one has
$-1$ as a repulsive fixed point, $1$ as attractive fixed point, and one can
easily check that every orbit besides $\left\{  1\right\}  $ and $\left\{
-1\right\}  $ for $\varphi$ is infinite, discrete and supported on a circle
passing by the points $-1$ and $1$. If $\varphi$ is elliptic, then up to
topological equivalence $\varphi$ is simply a rotation around the origin. If
$\varphi$ is parabolic, it has a unique fixed point $1$ (up to conformal
conjugacy) which is both attractive and repulsive, and its orbits besides
$\left\{  1\right\}  $ are infinite, discrete and supported on circles tangent
at $1$ to the vertical line of equation $x=1$.

\section{Hyperbolic Automorphisms}

Let $\varphi$ be a hyperbolic automorphism of $\mathbb{D}$ which fixes $1$ and
$-1$. By Theorem (\ref{topoconj}), any hyperbolic automorphism $\psi$ of
$\mathbb{D}$ is topologically conjugated to $\varphi$. Therefore, the
C*-crossed-products $C(\mathbb{D})\rtimes_{\varphi}\mathbb{Z}$ and
$C(\mathbb{D})\rtimes_{\psi}\mathbb{Z}$ are *-isomorphic. This section
describes the crossed-product $C(\mathbb{D})\rtimes_{\varphi}\mathbb{Z}$.

The orbit space of $\varphi$ restricted to $\mathbb{D}\backslash\left\{
-1,1\right\}  $ is denoted $\mathcal{C}$ and, topologically, is the cylinder
$[0,1]\times\mathbb{T}$. There is a one-to-one continuous map $\mathcal{D}%
_{\varphi}\rightarrow\mathcal{C}$, whose inverse is not continuous (note that
$\mathcal{C}$ is compact while $\mathcal{D}_{\varphi}$ is not).

\subsection{Irreducible representations}

We start with a full description of the irreducible *-representations of the
dynamical system $\left(  C(\mathbb{D)},\varphi,\mathbb{Z}\right)  $. We first
establish a result on the representations of general discrete dynamics
$\left(  C(\mathbb{D}),\psi,\mathbb{Z}\right)  $. We will denote the orbit for
$\psi$ of any $x\in\mathbb{D}$ by $\mathcal{O}_{\psi}(x)=\left\{  \psi
^{z}(x):z\in\mathbb{Z}\right\}  $.

\begin{lemma}
\label{spectral}Let $\psi$ be any homeomorphism of the closed unit disk
$\mathbb{D}$. Let $\pi$ be an irreducible representation of $C(\mathbb{D}%
)\rtimes_{\psi}\mathbb{Z}$ on some Hilbert space $\mathcal{H}$. Then
$\sigma(\pi(A_{\psi}))$ is a $\psi$-invariant compact subset of $\mathbb{D}$.
Moreover, for any $\psi$-invariant Borel subset $E$ of $\sigma(\pi(A_{\psi}%
))$, the spectral projection of $\pi(A_{\psi})$ associated to $E$ is either
$0$ or $1$. Hence, if $E$ is a nonempty relatively open invariant subset of
$\sigma(A_{\psi})$ then $\sigma(A_{\psi})=\overline{E}$.
\end{lemma}

\begin{proof}
Let $A_{\pi}=\pi(A_{\psi})$ and $U_{\pi}=\pi(U_{\psi})$. We have $U_{\pi
}^{\ast}A_{\pi}U_{\pi}=\psi(A_{\pi})$, so by the spectral mapping theorem
\cite[Theorem 8.11, p.289]{Conway90}, $\sigma(A_{\pi})$ is a $\psi$-invariant
set. Moreover, if we denote by $\chi_{E}$ the characteristic function of any
Borel subset $E$ of $\sigma(A_{\pi})$, then by the spectral theorem we have:%
\[
U_{\pi}^{\ast}\chi_{E}(A_{\pi})U_{\pi}=\chi_{E}\left(  U_{\pi}^{\ast}A_{\pi
}U\right)  =\chi_{\psi^{-1}(E)}(A_{\pi})\text{.}%
\]
In particular, assume $E$ is $\psi$-invariant. Then the spectral projection
$\chi_{E}(A_{\pi})$ of $A_{\pi}$ commutes with $U_{\pi}$ (and with $A_{\pi}$
by definition), hence since $\pi$ is irreducible we deduce that $\chi_{E}%
\in\left\{  0,1\right\}  $. If $E$ is relatively open in $\sigma(A_{\pi})$ and
nonempty, then $\chi_{E}(A_{\pi})>0$, so for any nonempty relatively open
$\psi$-invariant subset of $\sigma(A_{\pi})$ we have $\chi_{E}(A_{\pi})=1$. In
this case, $\chi_{\overline{E}}(A_{\pi})=1$, so $\chi_{\sigma(A_{\pi
})\backslash\overline{E}}(A_{\pi})=0$. Yet $\sigma(A_{\pi})\backslash
\overline{E}$ is relatively open in $\sigma(A_{\pi})$, so by the spectral
theorem $\sigma(A_{\pi})\backslash\overline{E}=\emptyset$.
\end{proof}

\begin{theorem}
\label{orbitspectrum}Let $\psi$ be an homeomorphism of the closed unit disk
$\mathbb{D}$. Let $\pi$ be an irreducible *-representation of $C(\mathbb{D}%
)\rtimes_{\psi}\mathbb{Z}$. Then there exists $x\in\mathbb{D}$ such that the
spectrum $\sigma\left(  \pi(A_{\psi})\right)  $ of $\pi(A_{\psi})$ is
$\overline{\mathcal{O}_{\psi}(x)}$.
\end{theorem}

\begin{proof}
Let $A_{\pi}=\pi(A_{\psi})$. For any subset $E$ of $\sigma(A_{\pi})$, we
denote the set $%
{\displaystyle\bigcup\limits_{n\in\mathbb{Z}}}
\psi^{n}\left(  E\right)  $ by $\mathcal{O}_{\psi}(E)$.

Let $\left(  \Omega_{n}\right)  _{n\in\mathbb{N}}$ be a countable basis of
relatively open nonempty subsets of $\sigma\left(  A_{\pi}\right)  $ for the
topology of the second countable compact space $\sigma(A_{\pi})$. Let $\Theta$
be a relatively open, nonempty subset of $\sigma(A_{\pi})$ and let
$n\in\mathbb{N}$. By Lemma (\ref{spectral}), we have $\sigma(A_{\pi
})=\overline{\mathcal{O}_{\psi}(\Omega_{n})}$, so $\Theta\cap\mathcal{O}%
_{\psi}(\Omega_{n})\not =\emptyset$. Now since $\sigma(A_{\pi})$ is a compact
Hausdorff space, there exists a nonempty open set $V$ such that $\overline
{V}\subseteq\Theta$, and we can choose $V\subseteq\mathcal{O}_{\psi}%
(\Omega_{n})$ for any choice of $n\in\mathbb{N}$ by Lemma (\ref{spectral}). We
can therefore construct by induction a sequence $\left(  V_{n}\right)
_{n\in\mathbb{N}}$ of relatively open nonempty subsets of $\sigma(A_{\pi})$
with $V_{0}=\Omega_{0}$ and such that for all $n\in\mathbb{N}$ we have
$\overline{V_{n+1}}\subseteq V_{n}$ and $V_{n}\subseteq\mathcal{O}_{\psi
}(\Omega_{n})$. Therefore, as $\sigma(A_{\pi})$ is compact, $%
{\displaystyle\bigcap\limits_{n\in\mathbb{N}}}
V_{n}=%
{\displaystyle\bigcap\limits_{n\in\mathbb{N}}}
\overline{V_{n}}$ is nonempty and by construction $%
{\displaystyle\bigcap\limits_{n\in\mathbb{N}}}
V_{n}\subseteq%
{\displaystyle\bigcap\limits_{n\in\mathbb{N}}}
\mathcal{O}_{\psi}(\Omega_{n})$. Let $x\in%
{\displaystyle\bigcap\limits_{n\in\mathbb{N}}}
V_{n}$: then since $x\in%
{\displaystyle\bigcap\limits_{n\in\mathbb{N}}}
\mathcal{O}_{\psi}(\Omega_{n})$, we have $\mathcal{O}_{\psi}(x)\cap
\mathcal{O}_{\psi}(\Omega_{n})\not =\emptyset$ for all $n\in\mathbb{N}$.
Therefore $\mathcal{O}_{\psi}(x)\cap\Omega_{n}\not =\emptyset$. Thus, we have
$\mathcal{O}_{\psi}(x)\cap\Theta\not =\emptyset$ since $\left(  \Omega
_{n}\right)  _{n\in\mathbb{N}}$ is a basis for the topology of $\sigma(A_{\pi
})$. Hence $\overline{\mathcal{O}_{\psi}(x)}=\sigma(A_{\pi})$.
\end{proof}

A way to read Lemma (\ref{spectral}) is to say that $\sigma(\pi(A_{\varphi}))$
is topologically transitive for any irreducible representation $\pi$.
Theorem\ (\ref{orbitspectrum}) holds in greater generality, for any second
countable Hausdorff compact space instead of $\mathbb{D}$, with essentially
the same proof. We only consider the case of the disk for our paper.

\bigskip We now work with our hyperbolic automorphism $\varphi$ of
$\mathbb{D}$. To fully describe the image of $C(\mathbb{D})\rtimes_{\varphi
}\mathbb{Z}$ by its irreducible representations, we will find the following
lemmas useful. We start by introducing some terminology: we let $\overline
{\mathbb{Z}}=\mathbb{Z}\cup\left\{  -\infty,\infty\right\}  $, and we let
$\mathbb{Z}$ acts by translation on $\overline{\mathbb{Z}}$. More precisely,
let $\tau$ be the action of $\mathbb{Z}$ on $\overline{\mathbb{Z}}$ defined by
$\tau_{z}(z^{\prime})=z+z^{\prime}$, where $z+\infty=\infty$ and
$z-\infty=-\infty$ for all $z\in\mathbb{Z}$ and $z^{\prime}\in\overline
{\mathbb{Z}}$. We define:

\begin{definition}
Let $\mathbb{H}=l^{2}\left(  \mathbb{Z}\right)  $ and let $\left(
e_{n}\right)  _{n\in\mathbb{Z}}$ be the canonical Hilbert basis of
$\mathbb{H}$. We define:%
\begin{align*}
C_{0,0}\left(  \mathbb{Z}\right)   &  =\left\{  \left(  \lambda_{n}\right)
_{n\in\mathbb{Z}}:\lim_{n\rightarrow\infty}\lambda_{n}=\lim_{n\rightarrow
-\infty}\lambda_{n}=0\right\}  \text{,}\\
C(\overline{\mathbb{Z}})  &  =\left\{  \left(  \lambda_{n}\right)
_{n\in\mathbb{Z}}:\left(  \lambda_{n}\right)  _{n\in\mathbb{N}}\text{ and
}\left(  \lambda_{-n}\right)  _{n\in\mathbb{N}}\text{ converge}\right\}
\text{.}%
\end{align*}
For any $s=\left(  \lambda_{n}\right)  _{n\in\mathbb{Z}}\in C(\overline
{\mathbb{Z}})$, let $\operatorname*{Diag}(s)$ be the unique bounded linear
operator on $\mathbb{H}$ such that $\operatorname*{Diag}(s)(e_{n})=\lambda
_{n}e_{n}$ for all $n\in\mathbb{Z}$. We set $\mathfrak{D}=\operatorname*{Diag}%
\left(  C\left(  \overline{\mathbb{Z}}\right)  \right)  $ and $\mathfrak{D}%
_{0,0}=\operatorname*{Diag}\left(  C_{0,0}\left(  \mathbb{Z}\right)  \right)
$. \ Of course, $\mathfrak{D}$ and $\mathfrak{D}_{0,0}$ are C*-algebras, and
$\operatorname*{Diag}$ is a *-isomorphism between $C(\overline{\mathbb{Z}})$
and $\mathfrak{D}$ and between $\mathfrak{D}_{0,0}$ and $C_{0,0}\left(
\mathbb{Z}\right)  $.
\end{definition}

\begin{definition}
\label{BilateralShiftDef}We define $U$ on $\mathbb{H}=l^{2}\left(
\mathbb{Z}\right)  $ by $U\left(  \xi_{n}\right)  _{n\in\mathbb{N}}=\left(
\xi_{n-1}\right)  _{n\in\mathbb{N}}$ for all $\xi\in\mathbb{H}$.\ Namely, $U$
is the bilateral shift unitary operator on $\mathbb{H}$.
\end{definition}

\begin{definition}
We endow $\overline{\mathbb{Z}}$ with the topology of the Gelfan'd spectrum of
$C\left(  \overline{\mathbb{Z}}\right)  $ (this is of course the usual
topology on $\overline{\mathbb{Z}}$).
\end{definition}

The next two lemmas describe the crossed-product $C\left(  \overline
{\mathbb{Z}}\right)  \rtimes_{\tau}\mathbb{Z}$, which can be seen as the
crossed-product $C(\overline{\mathcal{O}_{\varphi}(x)})\rtimes_{\varphi
}\mathbb{Z}$ when $x\in\mathbb{D}\backslash\left\{  -1,1\right\}  $ and is
therefore a fundamental building block of $C(\mathbb{D})\rtimes_{\varphi
}\mathbb{Z}$. Let $\mathcal{T}$ be the Toeplitz C*-algebra $C^{\ast}(S)$ where
$S$ is the unilateral shift on $l^{2}\left(  \mathbb{N}\right)  $ defined by
$S\left(  x_{0},x_{1},\ldots\right)  =\left(  0,x_{0},x_{1},\ldots\right)  $
for all $\left(  x_{n}\right)  _{n\in\mathbb{N}}\in l^{2}\left(
\mathbb{N}\right)  $.

\begin{lemma}
\label{UATTK}Let $E=\left\{  U^{n}\mathfrak{D}_{0,0}:n\in\mathbb{Z}\right\}
$. The C*-algebra $C^{\ast}\left(  E\right)  $ is the C*-algebra of compact
operators $\mathcal{K}$ and the C*-algebra $C^{\ast}\left(  U,\mathfrak{D}
\right)  $ equals $\left[
\begin{array}
[c]{cc}%
\mathcal{T} & \mathcal{K}\\
\mathcal{K} & \mathcal{T}%
\end{array}
\right]  $ where $\mathbb{H}=\mathbb{H}_{-}\oplus\mathbb{H}_{+}$ with
$\mathbb{H}_{-}=\left\{  e_{n}:n<0\right\}  $, $\mathbb{H}_{+}=\left\{
e_{n}:n\geq0\right\}  $ and with the natural identification $\mathbb{H}_{+}=\mathbb{H}_{-}
=l^{2}\left(  \mathbb{N}\right)  $. Moreover, if $\chi$ is defined by the
exact sequence $0\rightarrow\mathcal{K}\rightarrow\mathcal{T}\overset{\chi
}{\rightarrow}C(\mathbb{T})\rightarrow0$ and $\gamma=\chi\oplus\chi$ act on
$\left[
\begin{array}
[c]{cc}%
\mathcal{T} & \mathcal{K}\\
\mathcal{K} & \mathcal{T}%
\end{array}
\right]  $, then $\gamma(\operatorname*{Diag}(f)U^{n})=f(-\infty)Z^{-n}\oplus
f(\infty)Z^{n}$ for all $f\in C(\overline{\mathbb{Z}})$, $n\in\mathbb{Z}$ and
where $Z:z\in\mathbb{T}\mapsto z$.

Thus, $A\in\left[
\begin{array}
[c]{cc}%
\mathcal{T} & \mathcal{K}\\
\mathcal{K} & \mathcal{T}%
\end{array}
\right]  $ if and only if there exists $f,g\in C(\mathbb{T})$ and a compact
operator $k$ on $\mathbb{H}$ such that for any $\xi_{+}\oplus\xi_{-}%
\in\mathbb{H}_{-}\oplus\mathbb{H}_{+}$ we have $A\xi=\left(  T_{f}\xi
_{+}+T_{g}\xi_{-}\right)  +k\left(  \xi_{+}+\xi_{-}\right)  $ where $T_{f}$
and $T_{g}$ are the Toeplitz operators of symbols $f$ and $g$.
\end{lemma}

\begin{proof}
By construction, $\mathfrak{D}_{0,0}\subset\mathcal{K}$, so $C^{\ast}\left(
E\right)  \subseteq\mathcal{K}$ as $\mathcal{K}$ is an ideal in $\mathcal{B}%
\left(  \mathbb{H}\right)  $. Let now $P$ be any projection in $\mathcal{B}%
\left(  \mathbb{H}\right)  $. Assume that $P$ commutes with $C^{\ast}\left(
E\right)  $. Then $P$ commutes with $\mathfrak{D}_{0,0}$ so $P$ is diagonal in
the basis $\left(  e_{n}\right)  _{n\in\mathbb{Z}}$. If $P\not \in \{0,1\}$
then there exists $k\in\mathbb{Z}$ such that $Pe_{k}=0$ but $P_{e_{k-1}%
}=e_{k-1}$. Let $Q$ be the projection on $\operatorname*{span}\left\{
e_{n}:\left\vert n\right\vert \leq k\right\}  \in\mathfrak{D}_{0,0}$. Now, by
assumption $UQP=PUQ$ yet $UQPe_{k-1}=e_{k}$ but $PUe_{k-1}=0$. This is a
contradiction, so $P\in\{0,1\}$. Hence $C^{\ast}\left(  E\right)
=\mathcal{K}$ since it is irreducible \cite[Theorem 1.4.2, p. 18]{Arveson76}.

Write $U=\left[
\begin{array}
[c]{cc}%
U_{11} & 0\\
U_{21} & U_{22}
\end{array}
\right]  $ according to the decomposition $l^{2}\left(  \mathbb{Z}\right)
=\mathbb{H}_{-}\oplus\mathbb{H}_{+}$. Now, $U_{21}\in\mathcal{K}$ as it is the
rank one operator $U_{21}=\left\langle .,e_{-1}\right\rangle e_{0}$ (where
$\left\langle .,.\right\rangle $ is the inner product of $\mathbb{H}$). With
the obvious identifications of $\mathbb{H}_{-}$ and $\mathbb{H}_{+}$ with
$l^{2}\left(  \mathbb{N}\right)  $ it is immediate that $U_{11}=S^{\ast}$ and
$U_{22}=S$. Now, let $A$ be the projection on $\mathbb{H}_{+}$: by
construction, $A\in\mathfrak{D}$ and $AUA=U_{22}\in C^{\ast}(U,\mathfrak{D})$.
Similarly, $U_{11}\in C^{\ast}(U,\mathfrak{D})$, so $\mathcal{T}%
\oplus\mathcal{T}\subseteq C^{\ast}\left(  U,\mathfrak{D}\right)  $. Since
$C^{\ast}\left(  E\right)  \subset C^{\ast}\left(  U,\mathfrak{D}\right)  $ we
deduce that $\left[
\begin{array}
[c]{cc}%
\mathcal{T} & \mathcal{K}\\
\mathcal{K} & \mathcal{T}%
\end{array}
\right]  \subseteq C^{\ast}\left(  U,\mathfrak{D}\right)  $ (observe that
$\mathcal{K}$ is $\left[
\begin{array}
[c]{cc}%
\mathcal{K} & \mathcal{K}\\
\mathcal{K} & \mathcal{K}%
\end{array}
\right]  $ in our decomposition of $l^{2}\left(  \mathbb{Z}\right)  $). On the
other hand, by definition $C^{\ast}\left(  U,\mathfrak{D}\right)  $ is the
smallest C*-algebra containing $U$ and $\mathfrak{D}$. Since $U\in\left[
\begin{array}
[c]{cc}%
\mathcal{T} & \mathcal{K}\\
\mathcal{K} & \mathcal{T}%
\end{array}
\right]  $ and $\mathfrak{D}\subset\left[
\begin{array}
[c]{cc}%
\mathcal{T} & \mathcal{K}\\
\mathcal{K} & \mathcal{T}%
\end{array}
\right]  $ we deduce that $\left[
\begin{array}
[c]{cc}%
\mathcal{T} & \mathcal{K}\\
\mathcal{K} & \mathcal{T}%
\end{array}
\right]  =C^{\ast}\left(  U,\mathfrak{D}\right)  $. As a last remark, Let
$f\in C(\overline{\mathbb{Z}})$. Then $\operatorname*{Diag}%
(f)=(1-A)\operatorname*{Diag}(f)\oplus A\operatorname*{Diag}(f)$ and
$A\operatorname*{Diag}(f)=f(\infty)1+A\left(  \operatorname*{Diag}%
(f)-f(\infty)1\right)  $ where now $A\left(  \operatorname*{Diag}%
(f)-f(\infty)1\right)  \in\mathcal{K}$. Similarly, $(1-A)\operatorname*{Diag}%
(f)=f(-\infty)1+(1-A)(\operatorname*{Diag}(f)-f(-\infty)1)$ and
$(1-A)(\operatorname*{Diag}(f)-f(-\infty)1)\in\mathcal{K}$. Since
$\gamma(\mathcal{K})=\left\{  0\right\}  $ and $\gamma$ is a *-morphism with
$\gamma(U)=Z^{\ast}\oplus Z$, we deduce that $\gamma(\operatorname*{Diag}%
(f)U^{n})=f(-\infty)Z^{-n}\oplus f(\infty)Z^{n}$.
\end{proof}

\begin{lemma}
\label{ZZbarTTK}The C*-algebra $C\left(  \overline{\mathbb{Z}}\right)
\rtimes_{\tau}\mathbb{Z}$ is *-isomorphic to $\left[
\begin{array}
[c]{cc}%
\mathcal{T} & \mathcal{K}\\
\mathcal{K} & \mathcal{T}%
\end{array}
\right]  $.
\end{lemma}

\begin{proof}
Since $\mathbb{Z}$ is $\tau$--invariant, we have by \cite[Proposition
2.8]{Latremoliere05}:%
\[
0\longrightarrow C_{0,0}(\mathbb{Z})\rtimes_{\tau}\mathbb{Z}\longrightarrow
C\left(  \overline{\mathbb{Z}}\right)  \rtimes_{\tau}\mathbb{Z}\overset
{\gamma^{\prime}}{\longrightarrow}\mathbb{\ }C\left(  \left\{  -\infty
,\infty\right\}  \right)  \rtimes_{\tau}\mathbb{Z}\longrightarrow0\text{.}%
\]
Of course $C\left(  \left\{  -\infty,\infty\right\}  \right)  \rtimes_{\tau
}\mathbb{Z}=C(\mathbb{T})\oplus C(\mathbb{T})$ since both $\infty$ and
$-\infty$ are fixed point for $\tau$. Moreover, the C*-algebra $C_{0,0}%
(\mathbb{Z})\rtimes_{\tau}\mathbb{Z}\ $is *-isomorphic to $\mathcal{K}$
\cite[8.4.3, p. 141]{Fillmore96} and is the C*-algebra $C^{\ast}\left(
C_{0,0}(\mathbb{Z})V^{n}\right)  _{n\in\mathbb{Z}}$ where $V$ is the canonical
unitary in $C\left(  \overline{\mathbb{Z}}\right)  \rtimes_{\tau}\mathbb{Z}$.

Let us now define a map from $C\left(  \overline{\mathbb{Z}}\right)
\rtimes_{\tau}\mathbb{Z}$ onto $\left[
\begin{array}
[c]{cc}%
\mathcal{T} & \mathcal{K}\\
\mathcal{K} & \mathcal{T}%
\end{array}
\right]  $ as follows. Let $U$ be the bilateral shift in $\left[
\begin{array}
[c]{cc}%
\mathcal{T} & \mathcal{K}\\
\mathcal{K} & \mathcal{T}%
\end{array}
\right]  $. We set $\rho(V)=U$, and for $f\in C(\overline{\mathbb{Z}})$ we set
$\rho(f)=\operatorname*{Diag}(f)$. By universality of the crossed-product
$C\left(  \overline{\mathbb{Z}}\right)  \rtimes_{\tau}\mathbb{Z}$, we deduce
that $\rho$ extends into an epimorphism from $C\left(  \overline{\mathbb{Z}%
}\right)  \rtimes_{\tau}\mathbb{Z}$ onto $\left[
\begin{array}
[c]{cc}%
\mathcal{T} & \mathcal{K}\\
\mathcal{K} & \mathcal{T}%
\end{array}
\right]  $, since by construction $U^{\ast}\operatorname*{Diag}%
(f)U=\operatorname*{Diag}(g)$ where $g=f(.+1)$. Moreover, $\rho(fV^{n}%
)=\operatorname*{Diag}(f)U^{n}\in\mathcal{K}$ by construction, so $\rho\left(
C_{0,0}(\mathbb{Z})\rtimes_{\tau}\mathbb{Z}\right)  =C^{\ast}\left(  \left\{
U^{n}\mathfrak{D}_{0,0}:n\in\mathbb{Z}\right\}  \right)  =\mathcal{K}$ by
Lemma (\ref{UATTK}). Last, we observe that $\gamma^{\prime}%
(\operatorname*{Diag}(f)V^{n})=f(-\infty)Z^{\varepsilon_{1}n}\oplus
f(\infty)Z^{\varepsilon_{2}n}$ for all $f\in C(\overline{\mathbb{Z}})$,
$n\in\mathbb{Z}$, where $Z:z\in\mathbb{T}\mapsto z$ and with $\varepsilon
_{1},\varepsilon_{2}\in\left\{  -1,1\right\}  $. Up to a *-automorphism
of $C(\mathbb{T})\oplus C(\mathbb{T})$ we may as well choose $\varepsilon
_{1}=-1$ and $\varepsilon_{2}=1$. On the other hand, $\gamma(\rho
(\operatorname*{Diag}(f)V^{n}))=f(-\infty)Z^{-n}\oplus f(\infty)Z^{n}$ by
construction of $\rho$ and Lemma (\ref{UATTK}). Hence, we have the following
commuting diagram:%
\[%
\begin{array}
[c]{ccccccccc}%
0 & \rightarrow & C_{0}(\mathbb{Z})\rtimes_{\tau}\mathbb{Z} & \longrightarrow
& C\left(  \overline{\mathbb{Z}}\right)  \rtimes_{\tau}\mathbb{Z} &
\longrightarrow & C(\mathbb{T})\oplus C(\mathbb{T}) & \rightarrow & 0\\
&  & \downarrow\rho &  & \downarrow\rho &  & \parallel &  & \\
0 & \rightarrow & \mathcal{K} & \longrightarrow & \left[
\begin{array}
[c]{cc}%
\mathcal{T} & \mathcal{K}\\
\mathcal{K} & \mathcal{T}%
\end{array}
\right]  & \longrightarrow & C(\mathbb{T})\oplus C(\mathbb{T}) & \rightarrow &
0
\end{array}
\]
where the top and bottom lines are exact. Now, since $\mathcal{K}$ is simple,
$\rho$ in the first vertical arrow is injective, hence a *-isomorphism. The
last vertical arrow is just the identity of $C(\mathbb{T})\oplus
C(\mathbb{T})$. We deduce that $\rho$ in the second vertical arrow is a
*-isomorphism by the Short Five Lemma \cite[Lemma 3.1, p. 13]{Lane94}.
\end{proof}

We now define two families of *-representations of $C^{\ast}\left(
U_{\varphi},A_{\varphi}\right)  $, which we shall then prove in Theorem
(\ref{CovRepHyperbolic}) are the only ones up to unitary equivalence:

\begin{definition}
For any $\varepsilon\in\left\{  -1,1\right\}  $ and any $\theta\in\lbrack0,1)$
we define the representation $\pi_{\theta,\varepsilon}$ by $\pi_{\theta
,\varepsilon}(A_{\varphi})=\varepsilon$ and $\pi_{\theta,\varepsilon
}(U_{\varphi})=\exp\left(  2i\pi\theta\right)  $.

For any $x\in\mathbb{D}\backslash\left\{  -1,1\right\}  $, we define the
representation $\pi_{x}$ on $l^{2}\left(  \mathbb{Z}\right)  $ by letting
$\pi_{x}(A_{\varphi})=\operatorname*{Diag}\left(  \left(  \varphi
^{n}(x)\right)  _{n\in\mathbb{Z}}\right)  $ and $\pi_{x}(U_{\varphi})=U$ for
all $n\in\mathbb{Z}$.

We let $\operatorname*{IrrRep}=\left\{  \pi_{x},\pi_{\theta,\varepsilon}%
:x\in\mathbb{D}\backslash\left\{  -1,1\right\}  ,\varepsilon\in\left\{
-1,1\right\}  ,\theta\in\lbrack0,1)\right\}  $.
\end{definition}

We indeed check that for any $x\in\mathbb{D}\backslash\{-1,1\}$, $\theta
\in\lbrack0,1)$ and $\varepsilon\in\{-1,1\}$, if $\pi=\pi_{x}$ or $\pi
=\pi_{\theta,\varepsilon}\ $then we have $\pi(U_{\varphi}^{\ast}%
)\pi(A_{\varphi})\pi(U_{\varphi})=\varphi(\pi(A_{\varphi}))$, so $\pi$ extends
to a *-representation of $C^{\ast}(A_{\varphi},U_{\varphi})$ by universality. Now:

\begin{lemma}
All the representations in $\operatorname*{IrrRep}$ are irreducible. Let
$x,y\in\mathbb{D}\backslash\left\{  -1,1\right\}  ,$ $\varepsilon
,\varepsilon^{\prime}\in\{-1,1\}$ and $\theta,\theta^{\prime}\in\lbrack0,1)$.

We have $\pi_{x}\left(  C^{\ast}\left(  U_{\varphi},A_{\varphi}\right)
\right)  =\left[
\begin{array}
[c]{cc}%
\mathcal{T} & \mathcal{K}\\
\mathcal{K} & \mathcal{T}%
\end{array}
\right]  $, and moreover $\pi_{x}$ is unitarily equivalent to $\pi_{y}$ if and
only if $\mathcal{O}_{\varphi}(x)=\mathcal{O}_{\varphi}(y)$, and $\pi
_{\theta,\varepsilon}$ is unitarily equivalent to $\pi_{\varepsilon^{\prime
},\theta^{\prime}}$ if and only if $\left(  \varepsilon,\theta\right)
=\left(  \varepsilon^{\prime},\theta^{\prime}\right)  $.
\end{lemma}

\begin{proof}
For all representation $\pi$ the spectra of $\pi(A_{\varphi})$ and
$\pi(U_{\varphi})$ are invariant by unitary conjugation, it is also
straightforward that no two distinct elements of $\operatorname*{IrrRep}$ are
unitarily conjugated, and that if $\pi_{x}$ and $\pi_{y}$ are conjugated for
$x,y\in\mathbb{D}\backslash\{-1,1\}$ then $\overline{\mathcal{O}_{\varphi}%
(x)}=\sigma(\pi_{x}(A_{\varphi}))$ equals to $\sigma(\pi_{y}(A_{\varphi
}))=\overline{\mathcal{O}_{\varphi}(y)}$, and so $\mathcal{O}_{\varphi
}(x)=\overline{\mathcal{O}_{\varphi}(x)}\backslash\left\{  -1,1\right\}
=\overline{\mathcal{O}_{\varphi}(y)}\backslash\left\{  -1,1\right\}
=\mathcal{O}_{\varphi}(y)$. Conversely, if $\mathcal{O}_{\varphi
}(x)=\mathcal{O}_{\varphi}(y)$ then there exists $n\in\mathbb{Z}$ such that
$x=\varphi^{n}(y)$, so $U^{\ast n}\pi_{x}U^{n}=\pi_{y}$. The one-dimensional
representations $\pi_{\theta,\varepsilon}$, where $\theta\in\lbrack0,1)$ and 
$\varepsilon\in\{-1,1\}$, are obviously irreducible. If $x\in\mathcal{D}%
_{\varphi}$ then, by Lemma\ (\ref{UATTK}), we have $\pi_{x}\left(  C^{\ast
}\left(  U_{\varphi},A_{\varphi}\right)  \right)  =\left[
\begin{array}
[c]{cc}%
\mathcal{T} & \mathcal{K}\\
\mathcal{K} & \mathcal{T}%
\end{array}
\right]  $, since $C^{\ast}\left(  \pi_{x}(A_{\varphi})\right)  =\mathfrak{D}$
for all $x\in\mathbb{D}\backslash\left\{  -1,1\right\}  $. As in Lemma
(\ref{UATTK}), a projection $P$ of $\mathbb{H}$ commutes with $\left[
\begin{array}
[c]{cc}%
\mathcal{T} & \mathcal{K}\\
\mathcal{K} & \mathcal{T}%
\end{array}
\right]  $ if and only if $P\in\{0,1\}$ so $\pi_{x}$ is irreducible.
\end{proof}

\begin{theorem}
\label{CovRepHyperbolic}Let $\varphi$ be an hyperbolic conformal automorphism
of the closed unit disk $\mathbb{D}$. Let $\pi$ be an irreducible
representation for the dynamical system $\left(  C(\mathbb{D}),\varphi
,\mathbb{Z}\right)  $. Then $\pi$ is unitarily equivalent to one of the
representations in $\operatorname*{IrrRep}$.
\end{theorem}

\begin{proof}
Let $\pi$ be an irreducible representation of $\left(  C(\mathbb{D}%
),\varphi,\mathbb{Z}\right)  $ on some Hilbert space $\mathcal{H}$ with inner
product $\left\langle .,.\right\rangle $. Let $A_{\pi}=\pi(A_{\varphi})$ and
$U_{\pi}=\pi(U_{\varphi})$. Then $\sigma(A_{\pi})=\overline{\mathcal{O}%
_{\varphi}(x)}$ for some $x\in\mathbb{D}$ by Theorem (\ref{orbitspectrum}).
Now for our choice of automorphism $\varphi$, we observe that $x\in$
$\overline{\mathcal{O}_{\varphi}(x)}$ is always isolated in $\overline
{\mathcal{O}_{\varphi}(x)}$. Hence $x$ is an eigenvalue for $A_{\pi}$. Let
$\xi_{x}\in\mathcal{H}$ be a normalized eigenvector for $A_{\pi}$ associated
to $x$. Let $\mathcal{M}=\overline{\operatorname*{span}}\left\{  U_{\pi}%
^{n}\xi_{x}:n\in\mathbb{Z}\right\}  $. By construction, $\mathcal{M}$ is a
reducing subspace for both $U_{\pi}$ and $A_{\pi}$. Since $\pi$ is
irreducible, $\mathcal{M}\in\left\{  \mathcal{H},\left\{  0\right\}  \right\}
$, yet since $\xi_{x}\not =0$ we have $\mathcal{M}=\mathcal{H}$. Of course,
$A_{\pi}U_{\pi}^{n}\xi_{x}=U_{\pi}^{n}\varphi^{n}(A_{\pi})\xi_{x}=\varphi
^{n}(x)U_{\pi}^{n}\xi_{x}$; in other words $U_{\pi}^{n}\xi_{x}$ is a
normalized eigenvector for $A_{\pi}$ associated to the eigenvalue $\varphi
^{n}(x)$.

If $x\in\left\{  -1,1\right\}  $ then $\varphi^{n}(x)=x$ for all
$n\in\mathbb{Z}$. Hence $A_{\pi}=xI$ where $I$ is the identity of $\mathcal{H}$. 
Therefore, $C^\ast(A_\pi, U_\pi)$ is an Abelian C*-algebra, which is irreducible
by assumption on $\pi$. Hence it is one-dimensional and $\pi=\pi_{x,\theta}$ where $\theta \in [0,1)$ is given by $U_\pi \xi_x = \exp(2i\pi\theta) \xi_x$.

Now, if $x\not \in \{-1,1\}$ then let us define the linear operator
$W:l^{2}\left(  \mathbb{Z}\right)  \rightarrow\mathcal{M}$ by setting
$W(e_{n})=U_{\pi}^{n}\xi_{x}$ for all $n\in\mathbb{Z}$ and where $\left(
e_{n}\right)  _{n\in\mathbb{Z}}$ is the canonical basis of $l^{2}\left(
\mathbb{Z}\right)  $. Then $W^{\ast}\pi W=\pi_{x}$ by construction (since
$W^{\ast}\pi W\left(  A_{\varphi}\right)  (e_{n})=W^{\ast}A_{\pi}U_{\pi}%
^{n}\xi_{x}=\varphi^{n}(x)W^{\ast}\xi_{x}=\varphi^{n}(x)e_{n}$ and $W^{\ast
}\pi W(U_{\pi})(e_{n})=W^{\ast}U_{\pi}U_{\pi}^{n}\xi_{x}=W^{\ast}U_{\pi}%
^{n+1}\xi_{x}=e_{n+1}$). Now, $W$ is a unitary from $l^{2}\left(
\mathbb{Z}\right)  $ onto $\mathcal{M}$ since for all $n,m\in\mathbb{Z}$ we
have $\varphi^{n}(x)\not =\varphi^{m}(x)$ for $n\not =m$, hence $\left(
U_{\pi}^{n}\xi_{x}\right)  _{n\in\mathbb{Z}}$ is an orthonormal family. Hence,
$\pi$ is unitarily equivalent to $\pi_{x}$.
\end{proof}

\subsection{The Spectrum}

The spectrum $\widehat{C(\mathbb{D})\rtimes_{\varphi}\mathbb{Z}}$ of
$C(\mathbb{D})\rtimes_{\varphi}\mathbb{Z}$ can be identified with the set
$\operatorname*{IrrRep}$ by Theorem (\ref{CovRepHyperbolic}). We now describe
the topology on $\widehat{C(\mathbb{D})\rtimes_{\varphi}\mathbb{Z}}$. Since
all the representations of $\widehat{C(\mathbb{D})\rtimes_{\varphi}\mathbb{Z}%
}$ are of type\ I, the calculation of the spectrum could also use
\cite[Theorem 5.3]{Williams81}.

We recall:

\begin{proposition}
[{\cite[Theorem VIII.2.1, p. 222]{Davidson}}]\label{dualaction}Given any
crossed product $\mathcal{A}\rtimes_{f}\mathbb{Z}$ for some unital C*-algebra
$\mathcal{A}$ and $f$ any *-automorphism of $\mathcal{A}$, there exists a
unique action $\alpha$ of the Lie group $\mathbb{T}$ on $\mathcal{A}%
\rtimes_{f}\mathbb{Z}$ which satisfies $\alpha_{\lambda}(U_{f})=\lambda U_{f}$
and $\alpha_{\lambda}(h)=h$ for all $\lambda\in\mathbb{T}$, $h\in\mathcal{A}$
and where $U_{f}$ is the canonical unitary in $\mathcal{A}\rtimes
_{f}\mathbb{Z}$.
\end{proposition}

\begin{definition}
\label{Expectation}For any $a\in\mathcal{A}\rtimes_{f}\mathbb{Z}$ and
$n\in\mathbb{Z}$ we set \[\mathbb{E}_{n}^{\alpha}(a)=\int_{\mathbb{T}}
\alpha_{\lambda}(aU_{f}^{\ast n})d\mu(\lambda)\] where $\mu$ is the Haar
probability measure on $\mathbb{T}$.
\end{definition}

\begin{proposition}
[{\cite[Theorem VIII.2.2, p. 223]{Davidson}}]\label{Approx}The map
$\mathbb{E}_{n}^{\alpha}:\mathcal{A}\rtimes_{f}\mathbb{Z}\longrightarrow
\mathcal{A}$ is a surjection of norm 1 for all $n\in\mathbb{Z}$, and moreover
any $a\in\mathcal{A}\rtimes_{f}\mathbb{Z}$ is the limit of $\left(  \Sigma
_{k}^{\alpha}(a)=\sum_{n=-k}^{k}\left(  1-\frac{\left\vert n\right\vert }%
{k+1}\right)  \mathbb{E}_{n}^{\alpha}(a)U^{n}\right)  _{k\in\mathbb{N}}$.
\end{proposition}

We identify the Gelfan'd spectrum $\widehat{C(\mathbb{D})\rtimes_{\varphi
}\mathbb{Z}}$ of $C(\mathbb{D})\rtimes_{\varphi}\mathbb{Z}$ with the set
$\operatorname*{IrrRep}$ by Theorem (\ref{CovRepHyperbolic}):

\begin{definition}
Let $q:$ $\operatorname*{IrrRep}\longrightarrow\mathcal{C}\cup\left(
\mathbb{T}\times\left\{  -1,1\right\}  \right)  $, where $\mathcal{C}$ is the
orbit space of $\left(  C(\mathbb{D}\backslash\left\{  -1,1\right\}
),\varphi,\mathbb{Z}\right)  $, be defined by: $q(\pi_{x})=\mathcal{O}%
_{\varphi}(x)\in\mathcal{C}$ and $q\left(  \pi_{\theta,\varepsilon}\right)
=\left(  e^{2i\pi\theta},\varepsilon\right)  \in\mathbb{T}\times\left\{
-1,1\right\}  $
\end{definition}

\begin{definition}
We endow $\mathcal{C}\cup\left(  \mathbb{T}\times\left\{  -1,1\right\}
\right)  $ with the topology $\mathfrak{C}$ where a set $F$ is closed if and
only if it is either contained and closed in $\mathbb{T}\times\left\{
-1,1\right\}  $, or it is of the form $F^{\prime}\cup\left(  \mathbb{T}%
\times\left\{  -1,1\right\}  \right)  $ where $F^{\prime}$ is closed in
$\mathcal{C}$.
\end{definition}

\begin{theorem}
\label{spectrum}Let $\varphi$ be an hyperbolic automorphism of the closed unit
disk $\mathbb{D}$. The map $q$ is an homeomorphism of the spectrum
$\widehat{C(\mathbb{D})\rtimes_{\varphi}\mathbb{Z}}$ of the C*-crossed product
$C(\mathbb{D})\rtimes_{\varphi}\mathbb{Z}$ with $\left(  \mathcal{C}%
\cup\left(  \mathbb{T}\times\left\{  -1,1\right\}  \right)  ,\mathfrak{C}%
\right)  $.
\end{theorem}

\begin{proof}
We already know that $\widehat{C(\mathbb{D})\rtimes_{\varphi}\mathbb{Z}}$ is
identified, as a set, with $\operatorname*{IrrRep}$ by Theorem
(\ref{CovRepHyperbolic}).\ Since $\mathbb{D}\backslash\left\{  -1,1\right\}  $
is an open $\varphi$-invariant set, we have the exact sequence $$0\rightarrow
C(\mathbb{D}\backslash\left\{  -1,1\right\}  )\rtimes_{\varphi}\mathbb{Z}%
\rightarrow C(\mathbb{D})\rtimes_{\varphi}\mathbb{Z}\rightarrow C\left(
\left\{  -1,1\right\}  \right)  \rtimes_{\varphi}\mathbb{Z}\rightarrow0\textrm{.}$$ Let
$\operatorname*{IrrRep}_{1}$ be the subset of $\operatorname*{IrrRep}$
consisting of the one-dimensional representations $\pi_{\theta,\varepsilon}$
$\left(  \varepsilon\in\{-1,1\}\text{,}\theta\in\lbrack0,1)\right)  $, and let
$\operatorname*{IrrRep}_{\infty}$ be the subset of $\operatorname*{IrrRep}$
consisting of the infinite dimensional representations $\pi_{x}$
in $\operatorname*{IrrRep}$, where $x\in\mathbb{D}\backslash\{-1,1\}$.
Observe that $$q\left(  \operatorname*{IrrRep}_{1}\right)  =\mathbb{T}\times\left\{
-1,1\right\}  \;\;\textrm{and} \;\; q\left(  \operatorname*{IrrRep}_{\infty}\right)
=\mathcal{C}\textrm{.}$$

An irreducible *-representation $\pi$ is unitarily conjugated to some element
of $\operatorname*{IrrRep}_{1}$ if and only if $\pi$ vanishes on the ideal
$C(\mathbb{D}\backslash\left\{  -1,1\right\}  )\rtimes_{\varphi}\mathbb{Z}$.
Let $\pi_{\theta,\varepsilon}\in\operatorname*{IrrRep}_{1}$. Then the quotient
representation of $C(\left\{  -1,1\right\}  \times\mathbb{T})$ defined by
$\pi_{\theta,\varepsilon}$ is obviously $f\in C(\left\{  -1,1\right\}
\times\mathbb{T})\mapsto f\left(  \varepsilon,e^{2i\pi\theta}\right)  =f\circ
q\left(  \pi_{\theta,\varepsilon}\right)  $. By \cite[Proposition 3.2.1, p.
61]{Dixmier}, the map $q$ is an homeomorphism onto the spectrum of $C\left(
\left\{  -1,1\right\}  \right)  \rtimes_{\varphi}\mathbb{Z}$, and
$\operatorname*{IrrRep}_{1}$ is closed in $\widehat{C(\mathbb{D}%
)\rtimes_{\varphi}\mathbb{Z}}$. Hence, a subset $F$ of $\operatorname*{IrrRep}%
_{1}$ is closed in $\widehat{C(\mathbb{D})\rtimes_{\varphi}\mathbb{Z}}$ if and
only if $q(F)$ is closed.

Let $F$ be a nonempty closed set in $\widehat{C(\mathbb{D})\rtimes_{\varphi
}\mathbb{Z}}$. By \cite[Proposition 3.1.2, p. 60]{Dixmier}, there exists a
(nonempty) subset $\mathcal{F}$ of $C(\mathbb{D})\rtimes_{\varphi}\mathbb{Z}$
such that $F=\left\{  \pi\in\operatorname*{IrrRep}:\mathcal{F}\subseteq\ker
\pi\right\}  $. Assume that $\pi_{x}\in F$ for some $x\in\mathcal{D}_{\varphi
}$. The range of $\pi_{x}$ is $\left[
\begin{array}
[c]{cc}%
\mathcal{T} & \mathcal{K}\\
\mathcal{K} & \mathcal{T}%
\end{array}
\right]  $. We denote by $\chi$ the canonical surjection of $\mathcal{T}$ onto
$\mathcal{T}/\mathcal{K}=C(\mathbb{T})$. Let $P=1\oplus0$ and $Q=1-P$. For any
$\theta\in\lbrack0,1)$ and any $a\in C(\mathbb{D})\rtimes_{\varphi}\mathbb{Z}%
$, we set: $$\chi_{\theta}^{1}(a)=\chi(P\pi_{x}(a)P)(e^{2i\pi\theta}) \;\; \textrm{and} \;\;
\chi_{\theta}^{-1}(a)=\chi(Q\pi_{x}(a)Q)(e^{2i\pi\theta})\textrm{.}$$ Of course,
$\chi_{\theta}^{1}$ and $\chi_{\theta}^{-1}$ are one dimensional
*-representations of $C(\mathbb{D})\rtimes_{\varphi}\mathbb{Z}$ and so they
are unitarily equivalent to, respectively, $\pi_{1,\theta}$ and $\pi
_{-1,\theta}$. Hence if $\pi_{x}(a)=0$ then $\pi_{\theta,\varepsilon}(a)=0$
for all $\theta\in\lbrack0,1)$ and $\varepsilon\in\left\{  -1,1\right\}  $.
Hence, for all $a\in\mathcal{F}$ we have $\pi^{\prime}(a)=0$ for all $\pi$ for
all $\pi^{\prime}\in\operatorname*{IrrRep}_{1}$. Thus $F$ contains
$\operatorname*{IrrRep}_{1}$. Now, let $\pi\in F\cap\operatorname*{IrrRep}%
_{\infty}$: by Theorem (\ref{CovRepHyperbolic}), there exists $x\in
\mathcal{D}_{\varphi}$ such that $\pi=\pi_{x}$. We set: $$\mathcal{E}=\left\{
y\in\mathbb{D}\backslash\left\{  -1,1\right\}  :\exists n\in\mathbb{Z}%
\ \ \ \pi_{\varphi^{n}(y)}\in F\right\}  \textrm{.} $$ Denote by $\alpha$ the canonical
action of $\mathbb{T}$ on $C(\mathbb{D})\rtimes_{\varphi}\mathbb{Z}$ and
$\beta$ the canonical action of $\mathbb{T}$ on $C(\overline{\mathbb{Z}%
})\rtimes_{\varphi}\mathbb{Z}$. We then have, for all $\lambda\in\mathbb{T}$,
that $\pi_{x}\circ\alpha_{\lambda}=\beta_{\lambda}\circ\pi_{x}$ (as can be
checked on $A_{\varphi}$ and $U_{\varphi}$ and then extended to $C^{\ast
}(A_{\varphi},U_{\varphi})$ since $\pi_{x}\circ\alpha$ and $\beta\circ\pi_{x}$
are *-morphisms). For all $a\in\mathcal{F}$ we have $\pi_{y}(a)=0$, so
$\mathbb{E}_{n}^{\alpha}(a)(y)=\pi_{y}(\mathbb{E}_{n}^{\alpha}(a))=\int
_{\mathbb{T}}\beta_{\lambda}\circ\pi_{y}(aU_{\varphi}^{\ast n})d\mu
(\lambda)=0$ for all $y\in\mathcal{E}$ and for all $n\in\mathbb{Z}$. Thus, if
we set: $$\mathcal{F}^{\prime}=\left\{  a\in A:\forall n\in\mathbb{Z}\ \ \forall
y\in\mathcal{E}\ \ \ \mathbb{E}_{n}(a)(y)=0\right\}  $$ then $\left\{  \pi
\in\operatorname*{IrrRep}:\mathcal{F}^{\prime}\subseteq\ker\pi\right\}
\subseteq F$. Conversely, let $a\in C(\mathbb{D})\rtimes_{\varphi}\mathbb{Z}$
such that $\mathbb{E}_{n}^{\alpha}(a)(y)=0$ for all $y\in\mathcal{E}%
,n\in\mathbb{Z}$. Then $\Sigma_{k}^{\alpha}(a)(y)=\pi_{y}\left(  \Sigma
_{k}^{\alpha}(a)\right)  =0$ for all $y\in\mathcal{E},k\in\mathbb{N}$, and
thus $\pi_{y}(a)=0$ for all $y\in\mathcal{E}$ as $\pi_{y}$ continuous and
$\left(  \Sigma_{k}^{\alpha}(a)\right)  _{k\in\mathbb{N}}$ converges in norm
to $a$.\ Hence $\left\{  \pi\in\operatorname*{IrrRep}:\mathcal{F}^{\prime
}\subseteq\ker\pi\right\}  =F$. Now, let $\left(  \omega_{k}\right)
_{k\in\mathbb{N}}$ be a sequence in $\mathcal{C}$ such that $\omega_{k}\in
q(F\cap\operatorname*{IrrRep}_{\infty})$ for all $k\in\mathbb{N}$ and
converging to $\omega\in\mathcal{C}$. Let $a\in\mathcal{F}^{\prime}$. Then for
any $n\in\mathbb{Z}$ and any $y\in\omega$ we have $\mathbb{E}_{n}^{\alpha
}(a)(y)=0$ by continuity. Hence, $\omega\in q(F\cap\operatorname*{IrrRep}%
_{\infty}),$ so $q(F\cap\operatorname*{IrrRep}_{\infty})$ is a closed set in
$\mathcal{C}$. Therefore, $q(F)=\mathbb{T\cup}q(F\cap\operatorname*{IrrRep}%
_{\infty})$ is closed in $\mathfrak{C}$.

Conversely, let $F$ be a set in $\widehat{C(\mathbb{D})\rtimes
_{\varphi}\mathbb{Z}}$ such that $F\cap\operatorname*{IrrRep}_{\infty}
\not =\emptyset$ and $q(F)$ is closed. Then $T=q(F\cap\operatorname*{IrrRep}
_{\infty})$ is closed in $\mathcal{C}$. Let: $$\mathcal{M}=\left\{  a\in
C(\mathbb{D})\rtimes_{\varphi}\mathbb{Z}:\forall n\in\mathbb{Z}
\ \ \ \mathbb{E}_{n}(a)(T)=\left\{  0\right\}  \right\} \textrm{.} $$ Let $\pi
\in\operatorname*{IrrRep}_{\infty}$ such that $\pi(\mathcal{M})=\left\{
0\right\}  $ and let $\mathcal{O}=q(\pi)$. Since $\mathbb{E}_{0}^{\alpha}$ is
a surjection, for any $f\in C(\mathcal{C})$ such that $f(T)=\left\{
0\right\}  $ there exists $a\in\mathcal{M}$ such that $\mathbb{E}_{0}^{\alpha
}(a)=f$. Therefore, for all $f\in C\left(  \mathcal{C}\right)  $ such that
$f(T)=\left\{  0\right\}  $ we have $f(\mathcal{O})=0$, so $\mathcal{O}%
\in\overline{T}=T$. Hence $\pi\in F$. Hence, $F$ is closed in $\widehat
{C(\mathbb{D})\rtimes_{\varphi}\mathbb{Z}}$.
\end{proof}

A corollary of Theorem (\ref{spectrum}) is that the spectrum of the ideal
$\mathcal{I}=C(\mathbb{D}\backslash\left\{  -1,1\right\}  )\rtimes_{\varphi
}\mathbb{Z}$ is the set $\operatorname*{IrrRep}_{\infty}$ of infinite
dimensional representations in $\operatorname*{IrrRep}$, with the topology of
$\mathcal{C}$. Since all the representations in $\operatorname*{IrrRep}%
_{\infty}$ acts on the same Hilbert space $l^{2}\left(  \mathbb{Z}\right)  $,
we deduce moreover that $\mathcal{I}$ is a homogenous C*-algebra. Since
$\mathcal{I}$ has continuous trace, and since $H^{3}\left(  \mathcal{C}%
,\mathbb{Z}\right)  =\left\{  0\right\}  $, we deduce by \cite[Corollary
10.9.6, p. 219]{Dixmier} that $\mathcal{I}=C(\mathcal{C},\mathcal{K})$.
Moreover, it is easy to generalize the part of the proof of Theorem
(\ref{spectrum}) concerning $\mathcal{I}$ to any C*-crossed product
$C(X)\rtimes_{f}\mathbb{Z}$ where the action of $\mathbb{Z}$ on $X$ is proper
and free (so that $\mathcal{O}_{\varphi}(x)$ is homeomorphic to $\mathbb{Z}$
as a subset of $X$, for all $x\in X$): thus the spectrum of such
crossed-products is simply the orbit space of the action with its natural
quotient topology. We will provide a self-contained proof that $\mathcal{I}%
=C(\mathcal{C},\mathcal{K})$ in Proposition (\ref{idealsimpleform}), which
uses a similar but simplified method compared to \cite[Chapter 10]{Dixmier}.

\begin{remark}
Let $x\in\mathbb{D}\backslash\{-1,1\}$. Then by \cite[Proposition 3.2.1, p.
61]{Dixmier},we deduce from the exact sequence: $$0\rightarrow\ker\pi
_{x}\rightarrow C(\mathbb{D})\rtimes_{\varphi}\mathbb{Z}\overset{\pi_{x}%
}{\mathbb{\rightarrow}}\left[
\begin{array}
[c]{cc}
\mathcal{T} & \mathcal{K}\\
\mathcal{K} & \mathcal{T}%
\end{array}
\right]  \rightarrow0$$ that the spectrum of $\left[
\begin{array}
[c]{cc}
\mathcal{T} & \mathcal{K}\\
\mathcal{K} & \mathcal{T}%
\end{array}
\right]  $ is $\operatorname*{IrrRep}_{1}\cup\left\{  \pi_{x}\right\}  $,
topologized as a subset of $\left(  \operatorname*{IrrRep},\mathfrak{C}%
\right)  $. This result can be proven directly by using the exact sequence
$0\rightarrow\mathcal{K}\rightarrow\left[
\begin{array}
[c]{cc}%
\mathcal{T} & \mathcal{K}\\
\mathcal{K} & \mathcal{T}%
\end{array}
\right]  \rightarrow C(\mathbb{T})\oplus C(\mathbb{T})\rightarrow0$ and a
method similar yet simpler than the one used in Theorem (\ref{spectrum}).
\end{remark}

\subsection{The Crossed-Product}

The following C*-algebra is our candidate for the full crossed-product
C*-algebra $C(\mathbb{D})\rtimes_{\varphi}\mathbb{Z}$:

\begin{definition}
\label{A}We denote by $\mathcal{A}$ the C*-algebra:
\[
\left\{
\begin{array}
[c]{c}%
f\in C\left(  \mathbb{D}\backslash\left\{  -1,1\right\}  ,\left[
\begin{array}
[c]{cc}%
\mathcal{T} & \mathcal{K}\\
\mathcal{K} & \mathcal{T}%
\end{array}
\right]  \right)  \text{ such that:}\\
U^{\ast}fU=f\circ\varphi\text{ and }\forall t\in\mathbb{D}\backslash\left\{
-1,1\right\}  \ \ \ f(t)-f(0)\in\mathcal{K}%
\end{array}
\right\}
\]
where $U$ is the bilateral unitary shift defined by $U\left(  \zeta
_{n}\right)  _{n\in\mathbb{Z}}=\left(  \zeta_{n-1}\right)  _{n\in\mathbb{Z}}$
for all $\left(  \zeta_{n}\right)  _{n\in\mathbb{Z}}\in\mathbb{H}=l^{2}\left(
\mathbb{Z}\right)  $.
\end{definition}

The norm on $\mathcal{A}$ is $\left\Vert .\right\Vert _{\mathcal{A}}%
=\sup_{t\in\mathbb{D}\backslash\{-1,1\}}\left\Vert .\right\Vert _{C(\overline
{\mathbb{Z}})\rtimes_{\tau}\mathbb{Z}}$, which is indeed a norm (and then
easily checked to be a C*-norm)\ thanks to the condition $U^{\ast}%
fU=f\circ\varphi$, and is well-defined since $C(\overline{\mathbb{Z}}%
)\rtimes_{\tau}\mathbb{Z}=\left[
\begin{array}
[c]{cc}%
\mathcal{T} & \mathcal{K}\\
\mathcal{K} & \mathcal{T}%
\end{array}
\right]  $ by Lemma (\ref{ZZbarTTK}).

\begin{remark}
The C*-algebra $\mathcal{A}$ could alternatively be defined as the set of
$f\in C\left(  \mathbb{D}\backslash\left\{  -1,1\right\}  ,\left[
\begin{array}
[c]{cc}%
\mathcal{T} & \mathcal{K}\\
\mathcal{K} & \mathcal{T}%
\end{array}
\right]  \right)  $ such that $U^{\ast}fU=f\circ\varphi$ and $$\left(
\chi\oplus\chi\right)  (f(t))=\left(  \chi\oplus\chi\right)  (f(0))$$ for all
$t\in\mathbb{D}\backslash\left\{  -1,1\right\}  $ where $\chi:\mathcal{T}
\rightarrow C(\mathbb{T})$ is the canonical projection of the Toeplitz algebra
on $C(\mathbb{T})$ (i.e. the symbol map). This second condition reflects the fact
that $\varphi$ admits two and only two distinct fixed points on $\mathbb{D}$.
\end{remark}

\begin{lemma}
\label{Generators}Let $A^{\prime}(t)=\operatorname*{Diag}\left(  \left(
\varphi^{n}(t)\right)  _{n\in\mathbb{Z}}\right)  $ acting on $\mathbb{H}%
=l^{2}\left(  \mathbb{Z}\right)  $ for all $t\in\mathbb{D}\backslash\{-1,1\}$.
Then $A^{\prime}$ is a continuous function from $\mathbb{D}\backslash\{-1,1\}$
into the bounded operators on $\mathbb{H}$. Let $U^{\prime}:t\in
\mathbb{D}\backslash\left\{  -1,1\right\}  \mapsto U$ where $U\in\left[
\begin{array}
[c]{cc}%
\mathcal{T} & \mathcal{K}\\
\mathcal{K} & \mathcal{T}%
\end{array}
\right]  $, as defined in Definition (\ref{BilateralShiftDef}) and Definition
(\ref{A}), is the bilateral shift on $\mathbb{H}$. Then $C^{\ast}(A^{\prime
},U^{\prime})=\mathcal{A}$.
\end{lemma}

\begin{proof}
By definition $C^{\ast}(A^{\prime},U^{\prime})\subseteq\mathcal{A}$: indeed,
$A^{\prime}(t)-A^{\prime}(0)=0\oplus0$ and $U^{\prime}(t)-U^{\prime
}(0)=0\oplus0$ are independent of $t\in\mathbb{D}\backslash\left\{
-1,1\right\}  $. Moreover, $A\circ\varphi=U^{\ast}A^{\prime}U$ and $U^{\prime
}\circ\varphi=U^{\prime}=U^{\ast}U^{\prime}U$. Let us now prove the converse inclusion.

Let $\left(  e_{n}\right)  _{n\in\mathbb{Z}}$ be the canonical basis of
$l^{2}\left(  \mathbb{Z}\right)  $ and let $E_{n}$ be the orthogonal
projection in $l^{2}(\mathbb{Z})$ on $\mathbb{C}e_{n}$. Let $\mathcal{C}_{0}$
be the C*-algebra consisting of the operators $g(A^{\prime})$ for $g\in
C(\mathbb{D})$ such that $g(-1)=g(1)=0$, and let $\mathcal{A}_{0}=\left\{
f\in C\left(  \mathbb{D}\backslash\left\{  -1,1\right\}  ,\mathcal{K}\right)
:U^{\ast}fU=f\circ\varphi\right\}  $.

First, we observe that $\mathcal{C}_{0}$ is the C*-subalgebra $\mathcal{C}%
_{0}^{\prime}$ of $\mathcal{A}_{0}$ of the continuous functions on
$\mathbb{D}\backslash\left\{  -1,1\right\}  $ subject to the condition
$U^{\ast}fU=f\circ\varphi$ and such that $f(t)$ is diagonal in $\left(
e_{n}\right)  _{n\in\mathbb{Z}}$ for all $t\in$ $\mathbb{D}\backslash\left\{
-1,1\right\}  $. First, by construction $\mathcal{C}_{0}\subseteq
\mathcal{C}_{0}^{\prime}$. On the other hand, let $f\in\mathcal{C}_{0}%
^{\prime}$. For each $t\in\mathbb{D}\backslash\{-1,1\}$ we have $f(t)\in
C^{\ast}(A^{\prime}(t))$, so $f(t)=g_{t}(A^{\prime}(t))\ $for some continuous
function $g_{t}\in C(\overline{\mathcal{O}_{\varphi}\left(  t\right)  })$ such
that $g_{t}(1)=g_{t}(-1)=0$ (as $f(t)\ $is compact). Since $f$ is continuous,
we deduce that the function $g:t\in\mathbb{D}\backslash\left\{  -1,1\right\}
\mapsto g_{t}(t)$ is also continuous, and moreover it can be continuously
extended to $\mathbb{D}$ by setting $g(1)=g(-1)=0$. Now, since $U^{\ast
}A^{\prime}(t)U=A^{\prime}\circ\varphi(t)$ for all $t\in\mathbb{D}%
\backslash\left\{  -1,1\right\}  $, it follows that $g_{t}\left(
\varphi(t)\right)  =g_{\varphi(t)}(\varphi(t))$. Hence:%
\begin{align*}
g(A^{\prime}(t))  &  =g\left(  \operatorname*{Diag}\left(  \varphi
^{n}(t)\right)  _{n\in\mathbb{Z}}\right)  =\operatorname*{Diag}\left(  \left(
g\circ\varphi^{n}(t)\right)  _{n\in\mathbb{Z}}\right) \\
&  =\operatorname*{Diag}\left(  \left(  g_{\varphi^{n}(t)}\left(  \varphi
^{n}(t)\right)  \right)  _{n\in\mathbb{Z}}\right)  =\operatorname*{Diag}%
\left(  \left(  g_{t}\left(  \varphi^{n}(t)\right)  \right)  _{n\in\mathbb{Z}%
}\right) \\
&  =g_{t}(A^{\prime}(t))=f(t)\text{.}%
\end{align*}
Hence, $f\in\mathcal{C}_{0}$, so indeed $\mathcal{C}_{0}=\mathcal{C}%
_{0}^{\prime}$.

Now, let $f\in\mathcal{A}_{0}$. We write $f=f_{1}+f_{2}$ where $f_{1}%
\in\mathcal{C}_{0}$ such that $f_{1}(iy)=f(iy)$ for all $y\in\lbrack-1,1]$.
Thus, $f_{2}\in\mathcal{A}_{0}$ and moreover, $f_{2}\left(  iy\right)  =0$ for
all $y\in\lbrack-1,1]$. Hence, since $U^{\ast}f_{2}U=f_{2}\circ\varphi$, we
also have $f_{2}\left(  \varphi(iy)\right)  =0$. Thus, $f_{2}\in\left\{
g\in\mathcal{A}_{0}:g\left(  L\cup\varphi(L)\right)  =\left\{  0\right\}
\right\}  $ where $L=\left\{  iy:y\in\lbrack-1,1]\right\}  $. But it is now
immediate to check that $\left\{  g\in\mathcal{A}_{0}:g\left(  L\cup
\varphi(L)\right)  =\left\{  0\right\}  \right\}  $ is *-isomorphic to the
C*-subalgebra of operators $g\in C\left(  \overline{\mathcal{D}_{\varphi}%
}\right)  \otimes\mathcal{K}$ such that $g\left(  L\cup\varphi(L)\right)
=\left\{  0\right\}  $. Now, given any function $g\in C(\overline
{\mathcal{D}_{\varphi}})$, we define $g^{\prime}:t\in\mathbb{D}\backslash
\left\{  -1,1\right\}  \mapsto g([t])e_{n(t)}$ where $[t]$ and $n(t)$ are the
unique element in $\mathcal{D}_{\varphi}$ and $\mathbb{Z}$ respectively such
that $\varphi^{n(t)}\left(  [t]\right)  =t$. Now, it is easy to check that
$g^{\prime}$ is continuous on $\mathbb{D}\backslash\left\{  -1,1\right\}  $
(as the only problems are at the boundary, on which $g$ is zero everywhere),
and $U^{\ast}g^{\prime}U=g^{\prime}\circ\varphi$ as well as $g(t)\in
\mathcal{K}$ and is diagonal in $\left(  e_{n}\right)  _{n\in\mathbb{Z}}$ for
all $t\in\mathbb{D}\backslash\left\{  -1,1\right\}  $.\ Hence $g^{\prime}%
\in\mathcal{C}_{0}$. From this it is easy to check that $\left\{
g\in\mathcal{A}_{0}:g\left(  L\cup\varphi(L)\right)  =\left\{  0\right\}
\right\}  \subseteq C^{\ast}\left(  \mathcal{C}_{0},U^{\prime}\right)
=C^{\ast}(A^{\prime},U^{\prime})$. Therefore, $f_{2}\in C^{\ast}(A^{\prime
},U^{\prime})$. Hence $f\in C^{\ast}(A^{\prime},U^{\prime})$ and thus
$\mathcal{A}_{0}\subseteq C^{\ast}(A^{\prime},U^{\prime})$.

Now, we conclude our theorem as follows. Let $f\in\mathcal{A}$. There exists
$g\in C^{\ast}(A^{\prime},U^{\prime})$ such that $g(0)=f(0)$ by Lemma
(\ref{UATTK}), which proves that both $f(t)$ and $g(t)$ lie in $\left[
\begin{array}
[c]{cc}%
\mathcal{T} & \mathcal{K}\\
\mathcal{K} & \mathcal{T}%
\end{array}
\right]  $ for all $t\in\mathbb{D}\backslash\left\{  -1,1\right\}  $. Now,
$f(t)-g(t)\in\mathcal{K}$ so, by construction, $f-g\in\mathcal{A}_{0}$. Thus
$f-g\in C^{\ast}(A^{\prime},U^{\prime})$. This shows that $f=\left(
f-g\right)  +g\in C^{\ast}(A^{\prime},U^{\prime})$.

Hence as claimed, $\mathcal{A}=C^{\ast}(A^{\prime},U^{\prime})$.
\end{proof}

\begin{theorem}
\label{mainhyperbolic}Let $\varphi$ be a hyperbolic automorphism of the closed
unit disk $\mathbb{D}$. Then $C(\mathbb{D})\rtimes_{\varphi}\mathbb{Z}$ is
*-isomorphic to:%
\[
\mathcal{A}^{\prime}=\left\{
\begin{array}
[c]{c}%
f\in C\left(  \mathbb{R}\times\lbrack0,1],\left[
\begin{array}
[c]{cc}%
\mathcal{T} & \mathcal{K}\\
\mathcal{K} & \mathcal{T}%
\end{array}
\right]  \right)  \text{ such that:}\\
\forall x\in\mathbb{R}\ \forall y\in\lbrack0,1]\ \ \left\{
\begin{array}
[c]{c}%
f(x,y)-f(0,0)\in\mathcal{K}\\
U^{\ast}f(x,y)U=f(x+1,y)
\end{array}
\right.
\end{array}
\right\}  \text{.}%
\]

\end{theorem}

\begin{proof}
First we observe that $\sigma(A^{\prime})=\mathbb{D}$ and $U^{\prime\ast
}A^{\prime}U^{\prime}=\varphi(A^{\prime})$. Hence, there exists a unique
*-epimorphism $\theta$ from $C(\mathbb{D})\rtimes_{\varphi}\mathbb{Z}$ onto
$C^{\ast}(A^{\prime},U^{\prime})$ such that $\theta(A_{\varphi})=A^{\prime}$
and $\theta(U_{\varphi})=U^{\prime}$. Let $a\in C(\mathbb{D})\rtimes_{\varphi
}\mathbb{Z}$ such that $\theta(a)=0$. Let $\pi$ be an irreducible
*-representation of $C(\mathbb{D})\rtimes_{\varphi}\mathbb{Z}$. By Theorem
(\ref{CovRepHyperbolic}), up to unitary equivalence, $\pi\in
\operatorname*{IrrRep}$. Let $x\in\mathbb{D}$ such that $\sigma(\pi\left(
A_{\varphi}\right)  )=\overline{\mathcal{O}_{\varphi}(x)}$. If $x\in
\mathbb{D}\backslash\left\{  -1,1\right\}  $ then $\pi$ acts on $l^{2}\left(
\mathbb{Z}\right)  $ and extends: $\pi(A_{\varphi})f(z)=\varphi(x)^{z}f(z)$
and $\pi(U_{\varphi})f(z)=f(z-1).$ We set $\rho$ to be the faithful
*-representation of $\left[
\begin{array}
[c]{cc}%
\mathcal{T} & \mathcal{K}\\
\mathcal{K} & \mathcal{T}%
\end{array}
\right]  $ on $l^{2}\left(  \mathbb{Z}\right)  $ given in Lemma (\ref{UATTK}).
Let $\pi^{\prime}:f\in\mathcal{A}\mapsto\rho(f(x))$: by construction,
$\pi^{\prime}$ is an (irreducible) *-representation of $C^{\ast}(A^{\prime
},U^{\prime})$.

If $x=\{-1,1\}$ then $\pi$ acts on $\mathbb{C}$ by $\pi(U_{\varphi}%
)=\lambda\in\mathbb{T}$ and $\pi(A_{\varphi})=x$. If $x=1$ then we set
$\pi^{\prime}:f\in\mathcal{A}\mapsto(\left(  0\oplus\chi\right)
f(t))(\lambda)$ (for any $t\in\mathbb{D}\backslash\{-1,1\}$ since $\left(
0\oplus\chi\right)  (f(t))$ is independent of $t$ for $f\in\mathcal{A}$ by
construction). If $x=-1$ we set $\pi^{\prime}:f\in\mathcal{A}\mapsto\left(
\left(  \chi\oplus0\right)  f(t)\right)  \left(  \lambda\right)  $.

Whatever $x\in\mathbb{D}$ is, it is now easy to check that $\pi^{\prime}%
\circ\theta=\pi$, as this equality is valid on the total set of monomials in
$A_{\varphi},U_{\varphi}$. Therefore, $\pi(a)=\pi^{\prime}\circ\theta(a)=0$.
Since $\pi$ is arbitrary, we deduce $a=0$. Hence $\theta$ is also injective,
hence $\theta$ is a *-isomorphism from $C(\mathbb{D})\rtimes_{\varphi
}\mathbb{Z}$ onto:
\[
\mathcal{A}=\left\{
\begin{array}
[c]{c}%
f\in C\left(  \mathbb{D}\backslash\left\{  -1,1\right\}  ,\left[
\begin{array}
[c]{cc}%
\mathcal{T} & \mathcal{K}\\
\mathcal{K} & \mathcal{T}%
\end{array}
\right]  \right)  \text{ such that:}\\
U^{\ast}fU=f\circ\varphi\text{ and }\forall t\in\mathbb{D}\backslash\left\{
-1,1\right\}  \ \ \ f(t)-f(0)\in\mathcal{K}%
\end{array}
\right\}  \text{.}%
\]
Now, in the proof of Theorem (\ref{topoconj}), we constructed an homeomorphism
$\mu$ from $\mathcal{D}_{\varphi}$ onto $[0,1]^{2}$, which can be extended by
induction to an homeomorphism from $\mathbb{D}\backslash\left\{  -1,1\right\}
$ onto $\mathbb{R}\times\lbrack0,1]$ in such a way that\ $\mu^{-1}\circ
\varphi\circ\mu$ is simply the translation $(x,y)\in\mathbb{R}\times
\lbrack0,1]\mapsto(x+1,y)$. The map $\psi:f\in\mathcal{A}\mapsto f\circ\mu
\in\mathcal{A}^{\prime}$ is obviously a *-isomorphism.
\end{proof}

Since the action of $\mathbb{Z}$ on $\mathbb{D}\backslash\left\{
-1,1\right\}$ by an hyperbolic automorphism fixing $\{-1,1\}$ is proper and free, the ideal $C(\mathbb{D}\backslash\left\{
-1,1\right\}  )\rtimes_{\varphi}\mathbb{Z}$ of $C(\mathbb{D})\rtimes_{\varphi
}\mathbb{Z}$ is $C(\mathcal{C},\mathcal{K})$, as can be derived from our work
or the general theory of such crossed-products (which are known to be field of
elementary C*-algebras over the orbit space, thus classified by their Dixmier
-Douady class, which lives in $H^{3}\left(  \mathcal{C},\mathbb{Z}\right)
=\{0\}$, hence the field is trivial). We offer here the proof of this last
isomorphism as a corollary of our work:

\begin{proposition}
\label{idealsimpleform}Let $\varphi$ be a hyperbolic automorphism of
$\mathbb{D}$. Then $C_{-1,1}(\mathbb{D})\rtimes_{\varphi}\mathbb{Z}=C(\mathcal{C},\mathcal{K})$ where $\mathcal{C}%
=[0,1]\times\mathbb{T}$ and $C_{-1,1}(\mathbb{D})$ is the C*-algebra of continuous functions on $\mathbb{D}$ such that $f(-1)=f(1)=0$.
\end{proposition}

\begin{proof}
The C*-algebra $C_{-1,1}(\mathbb{D})\rtimes
_{\varphi}\mathbb{Z}$ is the kernel of the map $\eta$ by construction, so it
is isomorphic to $C\left(  \mathbb{D}\backslash\left\{  -1,1\right\}
,\mathcal{K}:U^{\ast}fU=f\circ\varphi\right)  $. There exists a selfadjoint
operator $H$ on $l^{2}\left(  \mathbb{Z}\right)  $ such that $U=\exp(iH)$. Let
$U_{t}=\exp(itH)$ for all $t\in\lbrack0,1]$. Let $h$ be an homeomorphism from
$\mathbb{D}\backslash\left\{  -1,1\right\}  \longrightarrow\mathbb{R}
\times\lbrack0,1]$ such that $h\varphi h^{-1}$ is just the translation by $1$
on $\mathbb{R}$. For any $f\in C(\mathcal{C},\mathcal{K})$ we define
$\widetilde{f}$ as follows: $\widetilde{f}(x)=U_{r}^{\ast}f\left(
e(r_{0}),s\right)  U_{r}$ where $h(r,s)=\omega$. Now, $\widetilde{f}\in
C\left(  \mathbb{D}\backslash\left\{  -1,1\right\}  ,\mathcal{K}:U^{\ast
}fU=f\circ\varphi\right)  $. The converse inclusion is shown similarly.
\end{proof}
\newpage
\subsection{K-Theory, Projections and Unitaries}

\subsubsection{K-theory for general actions on $\mathbb{D}$}

The whole $K$-theory of $\mathbb{D}$ can be read from the following six-terms
exact sequence \cite[Theorem 9.3.1, p. 67]{Blackadar98}:

\[
\begin{tabular}
[c]{lllll}
$K_{0}\left(  \overset{\circ}{\mathbb{D}}\right)  =\mathbb{Z}$ &
$\longrightarrow$ & $K_{0}\left(  \mathbb{D}\right)  =\mathbb{Z}$ &
$\longrightarrow$ & $K_{0}\left(  \mathbb{T}\right)  =\mathbb{Z}$\\
$\uparrow$ &  &  &  & $\downarrow$\\
$K_{1}\left(  \mathbb{T}\right)  =\mathbb{Z}$ & $\longleftarrow$ &
$K_{1}\left(  \mathbb{D}\right)  $ & $\longleftarrow$ & $K_{1}\left(
\overset{\circ}{\mathbb{D}}\right)  =0$
\end{tabular}
\]
which follows from the exact sequence $$0\longrightarrow C\left(
\overset{\circ}{\mathbb{D}}\right)  \longrightarrow C\left(  \mathbb{D}
\right)  \longrightarrow C\left(  \mathbb{T}\right)  \longrightarrow 0 \textrm{.}$$ We
conclude immediately that $K_{0}\left(  \mathbb{D}\right)  =\mathbb{Z}$ and
$K_{1}\left(  \mathbb{D}\right)  =0$. Yet, since $\mathbb{D}$ is a
contractible compact metric space, we immediately conclude that $K_{0}\left(
\mathbb{D}\right)  =\mathbb{Z}$ with generator $[1]$. From this, we can deduce
the $K$-theory of the crossed-product of $\varphi$ on $C\left(  \mathbb{D}%
\right)  $:

\begin{proposition}
We have $K_{0}\left(  C(\mathbb{D})\rtimes_{\varphi}\mathbb{Z}\right)
=\mathbb{Z}$, generated by the class of the identity, and $K_{1}\left(
C(\mathbb{D})\rtimes_{\varphi}\mathbb{Z}\right)  =\mathbb{Z}$ generated by the
canonical unitary of $C(\mathbb{D})\rtimes_{\varphi}\mathbb{Z}$, for any
homeomorphism $\varphi$ of $\mathbb{D}$.
\end{proposition}

\begin{proof}
The Pimsner-Voiculescu six-terms exact sequence
\cite{Pimsner80},\cite[Theorem 10.2.1, p. 73]{Blackadar98} for this action is:
\[
\begin{tabular}
[c]{lllll}
$K_{0}\left(  \mathbb{D}\right)  =\mathbb{Z}$ & $\overset{1-\varphi^{\ast}%
}{\longrightarrow}$ & $K_{0}\left(  \mathbb{D}\right)  =\mathbb{Z}$ &
$\longrightarrow$ & $K_{0}\left(  C(\mathbb{D})\rtimes_{\varphi}%
\mathbb{Z}\right)  $\\
$\uparrow$ &  &  &  & $\downarrow$\\
$K_{1}\left(  C(\mathbb{D})\rtimes_{\varphi}\mathbb{Z}\right)  $ &
$\longleftarrow$ & $K_{1}\left(  \mathbb{D}\right)  =0$ & $\overset
{1-\varphi^{\ast}}{\longleftarrow}$ & $K_{1}\left(  \mathbb{D}\right)  =0$%
\end{tabular}
\]
and deduce the result immediately since $\left(  1-\varphi^{\ast}\right)
[1]=0$. Note that the $K$-groups do not depend on the homeomorphism $\varphi$
of $\mathbb{D}$.
\end{proof}

\subsubsection{Another view of K-theory}

Even though the only nontrivial class in $K_{0}\left(  C(\mathbb{D}%
)\rtimes_{\varphi}\mathbb{Z}\right)  $ is always the one of the identity,
there can be many nontrivial projections in $C(\mathbb{D})\rtimes_{\varphi
}\mathbb{Z} $ whose $K_{0}$ classes vanish. This is illustrated by our
example, when $\varphi$ is an hyperbolic automorphism:

\begin{proposition}
Let $\varphi$ be a hyperbolic automorphism of $\mathbb{D}$. The projections of
$C(\mathbb{D})\rtimes_{\varphi}\mathbb{Z}$ of $K$-class $0$ are (in one-to-one
correspondence with) the projections in $C(\mathcal{C},\mathcal{K})$. Any
projection $P$ of $C(\mathbb{D})\rtimes_{\varphi}\mathbb{Z}$ is homotopic to
$Q$ or $1-Q$ for some compact--projection-valued continuous function $Q$ on
$\mathbb{D}\backslash\left\{  -1,1\right\}  $ such that $U^{\ast}%
QU=Q\circ\varphi$.

The unitaries of $C(\mathbb{D})\rtimes_{\varphi}\mathbb{Z}$ are all
homotopic\ to $U^{n}$ for some $n\in\mathbb{Z}$.
\end{proposition}

\begin{proof}
By Proposition (\ref{idealsimpleform}) we have the following exact sequence:%
\[
0\longrightarrow C(\mathcal{C},\mathcal{K})\longrightarrow C(\mathbb{D}%
)\rtimes_{\varphi}\mathbb{Z}\longrightarrow C\left(  \mathbb{T}\right)  \oplus
C\left(  \mathbb{T}\right)  \longrightarrow0
\]
which implies that the following six-terms exact sequence holds:%
\[%
\begin{tabular}
[c]{lllll}%
$K_{0}\left(  \mathcal{C},\mathcal{K}\right)  =\mathbb{Z}$ & $\longrightarrow$
& $K_{0}\left(  C(\mathbb{D})\rtimes_{\varphi}\mathbb{Z}\right)  =\mathbb{Z}$
& $\longrightarrow$ & $\mathbb{Z}^{2}$\\
$\uparrow$ &  &  &  & $\downarrow$\\
$\mathbb{Z}^{2}$ & $\longleftarrow$ & $K_{1}\left(  C(\mathbb{D}%
)\rtimes_{\varphi}\mathbb{Z}\right)  =\mathbb{Z}$ & $\longleftarrow$ &
$K_{1}\left(  \mathcal{C},\mathcal{K}\right)  =\mathbb{Z}$%
\end{tabular}
\ \text{.}%
\]
By inspecting each arrow, we conclude that the image of $K_{0}\left(
\mathcal{C},\mathcal{K}\right)  $ can only be $\left\{  0\right\}  $ for the
top row to be exact. Thus all compact projections in $C(\mathbb{D}%
)\rtimes_{\varphi}\mathbb{Z}$ have null $K_{0}$-class. Since obviously $1\in
C(\mathbb{D})\rtimes_{\varphi}\mathbb{Z}$ has for image $(1,1)\in
\mathbb{Z}^{2}$, we conclude that there is no projections in $C(\mathbb{D}%
)\rtimes_{\varphi}\mathbb{Z}$ which maps to either $0\oplus n$ or $n\oplus0$
for any $n\not =0$. Thus for any projection $P$ $\in$ $C(\mathbb{D}%
)\rtimes_{\varphi}\mathbb{Z}$ there exists a projection $K\in C(\mathbb{D}%
)\rtimes_{\varphi}\mathbb{Z}$ with $K(t)$ compact for all $t\in\mathbb{D}%
\backslash\left\{  -1,1\right\}  $ such that $P$ is homotopic (in the space of
projections) to $K$ or $1-K$.

Similarly, the bottom row exactness imposes that $K_{1}\left(  C(\mathbb{D}%
)\rtimes_{\varphi}\mathbb{Z}\right)  $ maps into a copy of $\mathbb{Z}$ in
$\mathbb{Z}^{2}$ and since $U$ has class $(1,-1)$ (since the image of $U$ in
$C(\mathbb{T})\oplus C(\mathbb{T})$ is by construction $z\mapsto
z\oplus\overline{z}$), we conclude that any unitary $W$ in $C(\mathbb{D}%
)\rtimes_{\varphi}\mathbb{Z}$ which $K_{1}$-class $n\in\mathbb{Z}$ is of the
form $V+K$ where $K\in C(\mathbb{D})\rtimes_{\varphi}\mathbb{Z}$ is
compact-valued (i.e. $\eta(K)=0$) and $V$ is homotopic (in the unitary
group)\ to $U^{n}$.
\end{proof}

In particular, one observes that if $\mathcal{A}^{\prime}$ is the C*-algebra:%
\[
\left\{  f\in C\left(  \mathcal{C},\left[
\begin{array}
[c]{cc}%
\mathcal{T} & \mathcal{K}\\
\mathcal{K} & \mathcal{T}%
\end{array}
\right]  \right)  :\left(  \chi\oplus\chi\right)  (f(t)-f(0))=0\right\}
\]
then $K_{0}\left(  \mathcal{A}^{\prime}\right)  =\mathbb{Z}^{2}$, so
Proposition (\ref{idealsimpleform}) does not extend to $C(\mathbb{D}%
)\rtimes_{\varphi}\mathbb{Z}$.

\section{Parabolic Automorphisms}

Informally, a parabolic automorphism of the closed unit disk $\mathbb{D}$ is a
limit case for the hyperbolic automorphisms we studied in the previous
section. However, the situation turns out to be easier in the parabolic case.
Let $\varphi$ be a parabolic automorphism of $\mathbb{D}$ which, by
Theorem\ (\ref{topoconj}), we can assume to fix $1$.

First, as in the hyperbolic case, there exists a domain $\mathcal{D}_{\varphi
}$ such that for all $x\in\mathbb{D}\backslash\left\{  1\right\}  $ there
exists a unique $y\in\mathcal{D}_{\varphi}$ such that $\mathcal{O}_{\varphi
}(x)=\mathcal{O}_{\varphi}(y)$. For instance, one can choose $\mathcal{D}%
_{\varphi}$ as the region between the diameter passing by $1$ and its image by
$\varphi$, not including this image. Second, observe that the orbit space of
$\varphi$ on $\mathbb{D}$ is the compact cone $\mathfrak{S}=\left\{  \left(
t,t\omega\right)  :t\in\lbrack0,1],\omega\in\mathbb{T}\right\}  $, which is
homeomorphic to $\mathbb{D}$.

\begin{definition}
Let $\varphi$ be a parabolic automorphism of $\mathbb{D}$ with fixed point $1$
(necessarily unique).

Let $\theta\in\lbrack0,1)$. We define the *-representation $\pi_{\theta}$ of
$C(\mathbb{D})\rtimes_{\varphi}\mathbb{Z}$ on $\mathbb{C}$ by setting:
$$\pi_{\theta}(A_{\varphi})=1 \;\; \textrm{and} \;\; \pi_{\theta}\left(  U_{\varphi}\right)
=\exp\left(  2i\pi\theta\right)  \textrm{.} $$

Let now $x\in\mathbb{D}\backslash\left\{  1\right\}  $. Then we define the
*-representation $\pi_{x}$ of $C(\mathbb{D})\rtimes_{\varphi}\mathbb{Z}$ on
$l^{2}\left(  \mathbb{Z},\right)  $ by setting: $$\pi_{x}(A_{\varphi})\left(
\left(  \xi_{n}\right)  _{n\in\mathbb{Z}}\right)  =\left(  \varphi^{n}%
(x)\xi_{n}\right)  _{n\in\mathbb{Z}} \;\; \textrm{and} \;\; \pi_{x}\left(  U_{\varphi}\right)
\left(  \left(  \xi_{n}\right)  _{n\in\mathbb{Z}}\right)  =\left(  \left(
\xi_{n-1}\right)  _{n\in\mathbb{Z}}\right)  $$ for all $\left(  \xi_{n}\right)
_{n\in\mathbb{Z}}\in l^{2}\left(  \mathbb{Z}\right) $.

We set $\operatorname*{IrrRep}=\left\{  \pi_{\theta},\pi_{x}:\theta\in
\lbrack0,1),x\in\mathcal{D}_{\varphi}\right\}  $.
\end{definition}

\begin{theorem}
\label{CovRepPara}Let $\varphi$ be a parabolic conformal automorphism of
$\mathbb{D}$. Then all the representations in $\operatorname*{IrrRep}$ are
irreducible *-representations, and no two representations in
$\operatorname*{IrrRep}$ are unitarily equivalent. Moreover, let $\pi$ be an
irreducible representation of $C(\mathbb{D})\rtimes_{\varphi}\mathbb{Z}$. Then
$\pi$ is unitarily equivalent to a representation in $\operatorname*{IrrRep}$.
\end{theorem}

\begin{proof}
It is straightforward to check that $\operatorname*{IrrRep}$ is indeed a set
of representations. If $\pi,\pi^{\prime}\in\operatorname*{IrrRep}$ and
$\pi\not =\pi^{\prime}$ then $\sigma(\pi(U_{\varphi}))\not =\sigma
(\pi(U_{\varphi}^{\prime}))$ or $\sigma(\pi\left(  A_{\varphi}\right)
)\not =\sigma\left(  \pi^{\prime}\left(  A_{\varphi}\right)  \right)  $ so
either way $\pi$ and $\pi^{\prime}$ are not unitarily conjugated.

\qquad Let $\pi$ be an irreducible representation of $C(\mathbb{D}%
)\rtimes_{\varphi}\mathbb{Z}$. By Theorem (\ref{orbitspectrum}) and denoting
$\pi(A_{\varphi})$ by $A_{\pi}$, we have $\sigma(A_{\pi})=\overline
{\mathcal{O}_{\varphi}(x)}$ for some $x\in\mathbb{D}$. The same proof as
Theorem (\ref{CovRepHyperbolic}) now applies: $x\in\overline{\mathcal{O}%
_{\varphi}(x)}$ is isolated in $\overline{\mathcal{O}_{\varphi}(x)}$ for
$\varphi$ parabolic, so there exists a vector $\xi_{x}$ of norm 1 such that
$A_{\pi}\xi_{x}=x\xi_{x}$. By irreducibility of $\pi$, it follows that
$\mathcal{H}=\overline{\operatorname*{span}}\left\{  U^{n}\xi_{n}%
:n\in\mathbb{Z}\right\}  $. If $x=1$ then $\mathcal{H}=\mathbb{C}$ and then
$\pi=\pi_{\theta}$ for $\theta\in\lbrack0,1)$. If $x\not =1$ then $\left(
U^{n}\xi_{n}\right)  _{n\in\mathbb{Z}}$ is orthonormal, and thus if we define
the operator $W:l^{2}\left(  \mathbb{Z}\right)  \rightarrow\mathcal{M}$ by
linearity such that $W\left(  e_{n}\right)  =U^{n}\xi_{x}$ for all
$n\in\mathbb{Z}$, where $\left(  e_{n}\right)  _{n\in\mathbb{Z}}$ is the
canonical basis for $l^{2}\left(  \mathbb{Z}\right)  $, then $W^{\ast}\pi
W=\pi_{x}$.
\end{proof}

We now identify the range of the infinite dimensional irreducible
representations of $C(\mathbb{D})\rtimes_{\varphi}\mathbb{Z}$:

\begin{lemma}
\label{shift}Let $x\in\mathbb{D}\backslash\left\{  1\right\}  $, and denote
$\pi_{x}(A_{\varphi})$ by $A_{x}$ and $\pi_{x}(U_{\varphi})$ by $U$. Then
$C^{\ast}\left(  U,A_{x}\right)  =C^{\ast}(U)+\mathcal{K}$ where $\mathcal{K}$
is the C*-algebra of compact operators on $l^{2}\left(  \mathbb{Z}\right)  $.
\end{lemma}

\begin{proof}
Let $x\in\mathbb{D}\backslash\left\{  1\right\}  $. Let $B\in\pi_{x}\left(
C(\mathbb{D})\rtimes_{\varphi}\mathbb{Z}\right)  =C^{\ast}\left(
A_{x},U\right)  $. Hence, the operator $B$ is the limit in norm of polynomials
in $A_{x}$ and $U$. Now, $A_{x}^{k}-1^{k}I$ is a compact operator, so $U^{\ast
n}\left(  A_{x}^{k}-1^{k}I\right)  U^{n}U^{m-n}$ is a compact operator. Hence
$C^{\ast}\left(  A_{x},U\right)  \subseteq\mathcal{K}+C^{\ast}(U)$.

Now, let $C$ be a finite-rank operator on $l^{2}\left(  \mathbb{Z}\right)  $:
then $C=\sum_{k=-n}^{n}C_{k}U^{k}$ where $C_{k}$ are finite rank diagonal
operators. Now, trivially $C_{k}=f_{k}(A_{x})$ for some continuous function
$f_{k}$ on $C(\overline{\mathcal{O}(x)})$ such that $f_{k}(1)=0$. Hence, $C\in
C^{\ast}(A_{x},U)$ and thus $\mathcal{K}+C^{\ast}(U)=C^{\ast}\left(
A_{x},U\right)  $.
\end{proof}

\begin{remark}
It is easy to prove that $C^{\ast}\left(  A_{x},U\right)  $ is the
C*-crossed-product $C(\mathbb{Z}\cup\left\{  \infty\right\}  )\rtimes_{\tau
}\mathbb{Z}$ where $\tau(z)=z+1$ and $\tau(\infty)=\infty$, and where
$\mathbb{Z}\cup\left\{  \infty\right\}  $ is the one-point compactification of
$\mathbb{Z}$.
\end{remark}

We now bundle all these irreducible representations together to form the
crossed-product C*-algebra:

\begin{theorem}
\label{mainpara}Let $\varphi$ be a parabolic conformal automorphism of the
closed unit disk $\mathbb{D}$ with fixed point $1\in\mathbb{T}$. Let
$\mathbb{H}=l^{2}\left(  \mathbb{Z}\right)  $ and let $U$ be the bilateral
shift $U\left(  \xi_{z}\right)  _{z\in\mathbb{Z}}=\left(  \xi_{z-1}\right)
_{z\in\mathbb{Z}}$ for all $\xi\in\mathbb{H}$. We also define, for all
$t\in\mathbb{D}\backslash\left\{  1\right\}  $, the operator $A_{x}(\xi
_{z})_{z\in\mathbb{Z}}=\left(  \varphi^{z}(x)\xi_{z}\right)  _{z\in\mathbb{Z}%
}$. Let $U^{\prime}:t\in\mathbb{D}\backslash\left\{  1\right\}  \mapsto U$ and
$A^{\prime}:t\in\mathbb{D}\backslash\left\{  1\right\}  \mapsto A_{t}$. Then
$C(\mathbb{D})\rtimes_{\varphi}\mathbb{Z}$ is *-isomorphic to the C*-algebra
$C^{\ast}(A^{\prime},U^{\prime})$, or equivalently to the C*-algebra:%
\[
\mathcal{A}=C^{\ast}(U^{\prime})+C_{1}(\mathbb{D},\mathcal{K})
\]
where $C_{1}(\mathbb{D},\mathcal{K})$ is the C*-algebra of $\mathcal{K}%
$-valued continuous functions on $\mathbb{D}$ vanishing at $1$.
\end{theorem}

\begin{proof}
Since $A^{\prime}=I+K=U^{\prime}U^{\prime\ast}+K$ where $K\in C_{1}%
(\mathbb{D},\mathcal{K})$, it is obvious that $C^{\ast}(A^{\prime},U^{\prime
})\subseteq C^{\ast}(U^{\prime})+C_{1}(\mathbb{D},\mathcal{K})$. We can
proceed as in Lemma\ (\ref{Generators}) to prove that $C^{\ast}(A^{\prime
},U^{\prime})=C^{\ast}(U^{\prime})+C_{1}(\mathbb{D},\mathcal{K})$: if $f\in
C^{\ast}(U^{\prime})+C_{1}(\mathbb{D},\mathcal{K})$ then $f=g+h$ where $g\in
C^{\ast}(U^{\prime})$ and $h\in C_{1}(\mathbb{D},\mathcal{K})$. Now
$C_{1}(\mathbb{D},\mathcal{K})=\left\{  h\in C(\mathbb{D}):h(1)=0\right\}
\otimes\mathcal{K}$. Yet, $C^{\ast}(A^{\prime})$ is the C*-algebra of
continuous functions on $\mathbb{D}$ vanishing at $1$ and valued in the
diagonal compact operators for the canonical basis of $l^{2}\left(
\mathbb{Z}\right)  $. Thus, $C_{1}(\mathbb{D},\mathcal{K})\subseteq C^{\ast
}(A^{\prime},U^{\prime})$ trivially, and thus $f\in C^{\ast}(A^{\prime
},U^{\prime})$ as desired.

Therefore $\mathcal{A}$ is generated by $U^{\prime}$ and $A^{\prime}$ such
that $U^{\prime\ast}A^{\prime}U^{\prime}=\varphi(A^{\prime})$ and
$\sigma(A^{\prime})=\mathbb{D}$, thus there exists a *-epimorphism
$\psi:C(\mathbb{D})\rtimes_{\varphi}\mathbb{Z}\twoheadrightarrow\mathcal{A}$
such that $\psi(A_{\varphi})=A^{\prime}$ and $\psi(U_{\varphi})=U^{\prime}$.
Now, let $\pi$ be an irreducible *-representation of $C(\mathbb{D}%
)\rtimes_{\varphi}\mathbb{Z}$. Up to unitarily equivalence, $\pi
\in\operatorname*{IrrRep}$. Let $a\in C(\mathbb{D})\rtimes_{\varphi}%
\mathbb{Z}$ such that $\psi(a)=0$. Assume first $\pi=\pi_{\theta}$ for some
$\theta\in\lbrack0,1)$. Define $\pi^{\prime}(U^{\prime})=e^{2i\pi\theta}$ and
$\pi^{\prime}(A^{\prime})=1$. Thus $\pi^{\prime}$ extends to a
*-representation of $\mathcal{A}$ and by construction, $\pi^{\prime}\circ
\psi=\pi$, so $\pi(a)=0$. Assume now $\pi=\pi_{x}$ for $x\in\mathbb{D}%
\backslash\left\{  1\right\}  $, as defined in Lemma (\ref{shift}). Then set
$\pi^{\prime}(A_{\varphi})=A_{x}$ and $\pi^{\prime}(U_{\varphi})=U$. Again, we
have $\pi(a)=\pi^{\prime}\circ\psi(a)=0$. Hence $\psi$ is injective, which
concludes our proof.
\end{proof}

We conclude this section with a description of $\widehat{C(\mathbb{D}%
)\rtimes_{\varphi}\mathbb{Z}}$. We work out the spectrum from Theorem
(\ref{mainpara}) rather than from pure representation theory considerations as
in Theorem (\ref{spectrum}). Let $\mathbb{D}$ be the orbit space of $\left(
C(\mathbb{D}),\varphi,\mathbb{Z}\right)  $ with its natural quotient topology.

\begin{theorem}
\label{paraspectrum}Let $\varphi$ be a parabolic automorphism of the closed
unit disk $\mathbb{D}$. We define the map $q:\operatorname*{IrrRep}%
\longmapsto\mathbb{T}\cup\mathbb{D}\backslash\{1\}$ by setting: $q\left(
\pi_{\theta}\right)  =\lambda$ and $q\left(  \pi_{x}\right)  =x$. We define on
$\mathbb{T}\cup\mathbb{D}\backslash\{1\}$ the topology $\mathfrak{C}$ where a
subset $F\subseteq\mathbb{T}\cup\mathbb{D}\backslash\{1\}$ is closed if and
only if it is either contained and closed in $\mathbb{T}$, or it is of the
form $F^{\prime}\cup\mathbb{T}$ where $F^{\prime}$ is closed in $\mathbb{D}%
\backslash\{1\}$. The spectrum $\widehat{C(\mathbb{D})\rtimes_{\varphi
}\mathbb{Z}}$ of $C(\mathbb{D})\rtimes_{\varphi}\mathbb{Z}$ is homeomorphic to
$\left(  \mathbb{T}\cup\mathbb{D}\backslash\{1\},\mathfrak{C}\right)  $.
\end{theorem}

\begin{proof}
By Theorem (\ref{mainpara}), we have $C(\mathbb{D})\rtimes_{\varphi
}\mathbb{Z=}C^{\ast}(U^{\prime})+C_{1}(\mathbb{D},\mathcal{K})$. For any $f\in
C(\mathbb{D})\rtimes_{\varphi}\mathbb{Z}$ we write $f=f_{U}+f_{\mathcal{K}}$
where $f_{U}\in C^{\ast}(U^{\prime})$ and $f_{\mathcal{K}}\in C_{1}%
(\mathbb{D},\mathcal{K})$ are uniquely determined by this equation. Let $\pi$
be an irreducible *-representation of $C^{\ast}(U^{\prime})+C_{1}%
(\mathbb{D},\mathcal{K})$. Then $\pi$ restricted to the ideal $C_{1}%
(\mathbb{D},\mathcal{K})$ is also irreducible, so it is either $0$ or the
evaluation map at some $x\in\mathbb{D}\backslash\left\{  1\right\}  $. In the
first case, then $\pi$ is really an irreducible representation of $C^{\ast
}(U)$, namely $\pi=\pi_{\theta}$ for some $\theta\in\lbrack0,1)$, and in the
second case $\pi=\pi_{x}$. Now, we have: $\pi_{\theta}(f)=f_{U}(e^{2i\pi
\theta})$ and $\pi_{x}(f)=f_{U}+f_{\mathcal{K}\lfloor\mathcal{O}(x)}$ where
$f_{\mathcal{K}\lfloor\mathcal{O}(x)}$ is the restriction of $f_{\mathcal{K}}$
to $\mathcal{O}(x)$.

Let $F$ be a subset of $\widehat{C(\mathbb{D})\rtimes_{\varphi}\mathbb{Z}}$.
Let $f\in\cap_{\pi\in F}\ker\pi$. Suppose first that there exists
$x\in\mathbb{D}\backslash\left\{  1\right\}  $ such that $\pi_{x}\in F$. Let
$y$ in the closure of $G=\left\{  \mathcal{O}(x)\in\mathbb{D}:\pi_{x}\in
F\right\}  $ for $\mathfrak{C}$. Then either $y\in\mathbb{T}$ or
$y\in\overline{G}\cap\mathcal{O}$. Now, $\pi_{x}(f)=f_{U}+f_{\mathcal{K}%
\lfloor\mathcal{O}(x)}=0$, so $f_{U}=0$ and thus $\pi_{y}(f)=0$ for all
$y\in\mathbb{T}$. On the other hand, since $f_{\mathcal{K}}$ is continuous,
$f_{\mathcal{K}\lfloor y}=0$ if $y$ is the limit in $\mathbb{D}$ of elements
in $G$. Thus, $F$ is closed and contains an infinite dimensional irreducible
representation if and only if $y\in\overline{G}$, hence if and only if $G$ is
closed in $\mathfrak{C}$. On the other hand if $F$ contains only one
dimensional representations, then $F=\left\{  \pi_{\theta}:e^{2i\pi\theta}\in
G\right\}  $ for some $G\subseteq\mathbb{T}$ is closed if and only if $G$ is
closed. This concludes the description of the topology on $\widehat
{C(\mathbb{D})\rtimes_{\varphi}\mathbb{Z}}$.
\end{proof}

\section{Elliptic Automorphisms}

The elliptic maps naturally give rise to slight generalization of the rotation
algebras. For any $\theta\in\lbrack0,1]$ we denote by $\mathcal{A}_{\theta}$
the universal C*-algebra generated by two unitaries $U$ and $V$ such that
$VU=e^{2i\pi\theta}UV$ \cite{Rieffel81}. We first observe that $C(\mathbb{D}%
)\rtimes_{\varphi}\mathbb{Z}$ is by definition the universal C*-algebra
generated by one unitaries $U_{\varphi}$ and one normal element $A_{\varphi}$
such that $\sigma(A_{\varphi})=\mathbb{D}$ and $U_{\varphi}^{\ast}A_{\varphi
}U_{\varphi}=e^{2i\pi\theta}A_{\varphi}$.

\begin{proposition}
\label{CovRepEll}Let $\theta\in\lbrack0,1)$, and for any $r\in\lbrack0,1]$ we
denote by $\Gamma_{r}$ the circle of center $0$ and radius $r$. We let
$\xi_{r}$ be defined as the quotient map:%
\[
C(\mathbb{D})\rtimes_{\varphi}\mathbb{Z}\overset{\xi_{r}}{\longrightarrow
}\left\{
\begin{array}
[c]{c}%
C(\Gamma_{r})\rtimes_{\varphi}\mathbb{Z=}\mathcal{A}_{\theta}\text{ if
}r>0\text{,}\\
C(\Gamma_{0})\rtimes_{\varphi}\mathbb{Z=}C(\mathbb{T})\text{ if }r=0\text{.}%
\end{array}
\right.
\]

Then $\pi$ is an irreducible representation of $C(\mathbb{D})\rtimes_{\varphi
}\mathbb{Z}$ if and only if either of the following holds:

\begin{itemize}
\item There exists $r\in]0,1]$ and $\rho$ an irreducible representation of
$\mathcal{A}_{\theta}$ such that $\pi=\rho\circ\xi_{r}$,

\item There exists $\omega\in\mathbb{T}$ such that $\pi=\nu_{\omega}\circ
\xi_{0}$, where $\nu_{\omega}(f)=f(\omega)$ for all $f\in C(\mathbb{T})$.
\end{itemize}

Consequently, two representations $\rho\circ\xi_{r}$ and $\rho^{\prime}%
\circ\xi_{r^{\prime}}$ are unitarily equivalent if and only if $r=r^{\prime}$
and $\rho$ and $\rho^{\prime}$ are unitarily equivalent representations of
$C(\Gamma_{r})\rtimes_{\varphi}\mathbb{Z}$.
\end{proposition}

\begin{proof}
Since $\varphi(\Gamma_{r})=\Gamma_{r}$ for all $r\in\lbrack0,1]$, the
following exact sequence holds:%
\[
0\longrightarrow C(\mathbb{D}\backslash\Gamma_{r})\rtimes_{\varphi}%
\mathbb{Z}\longrightarrow C(\mathbb{D})\rtimes_{\varphi}\mathbb{Z}\overset
{\xi_{r}}{\longrightarrow}C(\Gamma_{r})\rtimes\mathbb{Z}\longrightarrow0
\]
so $\xi_{r}$ if well-defined. Let $\pi$ be an irreducible *-representation of
$C(\mathbb{D})\rtimes_{\varphi}\mathbb{Z}$. Then by Theorem
(\ref{orbitspectrum}), there exists $r\in\lbrack0,1]$ such that $\sigma
(\pi(A_{\varphi}))\subseteq\Gamma_{r}$. If follows naturally that $\pi\left(
C(\mathbb{D}\backslash\Gamma_{r})\rtimes_{\varphi}\mathbb{Z}\right)  =\left\{
0\right\}  $ and thus there exists a *-representation $\rho$ of $C(\Gamma
_{r})\rtimes\mathbb{Z}$ such that $\pi=\rho\circ\xi_{r}$. Since $\pi$ is
irreducible, so is $\rho$. If $r>0$ then $C(\Gamma_{r})\rtimes\mathbb{Z}%
=\mathcal{A}_{\theta}$, while if $r=0$ then $C(\Gamma_{0})\rtimes
\mathbb{Z}=C(\mathbb{T})$, and this concludes the necessary condition. Now,
the conditions are trivially sufficient, and the unitary conjugacy
classification follows immediately.
\end{proof}

We now provide the description of the crossed-product C*-algebra for elliptic
automorphisms, starting by the following lemma:

\begin{lemma}
\label{fieldsofAtheta}Let $A^{\prime}:t\in\lbrack0,1]\mapsto tV\in
\mathcal{A}_{\theta}$ and $U^{\prime}:t\in\lbrack0,1]\mapsto U\in
\mathcal{A}_{\theta}$. Then $C^{\ast}(A^{\prime},U^{\prime})=\left\{  f\in
C\left(  [0,1],\mathcal{A}_{\theta}\right)  :f(0)\in C^{\ast}(U)\right\}  $.
\end{lemma}

\begin{proof}
Let $f\in C\left(  [0,1],\mathcal{A}_{\theta}\right)  $ such that $f(0)\in
C^{\ast}(U)$. We write $f=f_{1}+f_{2}$ where $f_{1}:t\in\lbrack0,1]\mapsto
f(0)\in C^{\ast}(U)$. Obviously, $f_{1}\in C^{\ast}(U^{\prime})$. Now,
$f_{2}\in C_{0}\left(  [0,1]\right)  \otimes\mathcal{A}_{\theta}$. Observe
that $t\in\lbrack0,1]\mapsto t^{2}=A^{\prime\ast}A^{\prime}(t)\in C^{\ast
}(A^{\prime},U^{\prime})$. By the Stone-Weierstrass theorem, we deduce that
$C^{\ast}(A^{\prime\ast}A^{\prime})=C_{0}\left(  [0,1]\right)  $. Let
$\varepsilon>0$ and $m\in\mathbb{Z}$, and let $\delta>0$ such that
$g(t)\leq\varepsilon$ for $t\in\lbrack0,\delta]$. Let $h_{\delta}\in
C_{0}\left(  [0,1]\right)  $ such that $h_{\delta}(t)=t^{-m}$ if $t\geq\delta
$, while $h_{\delta}(t)\leq\frac{1}{t^{m}}$ if $t\in\lbrack0,\delta]$ with
$h_{\delta}(0)=0$. Then $\left\Vert g(t)V^{m}-g(t)h_{\delta}(t)t^{m}%
V^{m}\right\Vert =0$ if $t\geq\delta$ and if $t\in\lbrack0,\delta]$ then
$\left\Vert g(t)V^{m}-g(t)h_{\delta}(t)t^{m}V^{m}\right\Vert \leq\left\vert
g(t)\right\vert \left\vert 1-h_{\delta}(t)t^{m}\right\vert \leq\varepsilon$.
Since $h_{\delta}(t)\in C^{\ast}(A^{\prime},U^{\prime})$ and $t\mapsto
t^{m}V^{m}=A^{\prime m}$ we deduce that $gV^{m}\in C^{\ast}(A^{\prime
},U^{\prime})$ since $C^{\ast}(A^{\prime},U^{\prime})$ is norm closed. Thus we
deduce immediately that for all$k\in\mathbb{Z}$ we have $gV^{m}U^{k}\in
C^{\ast}(A^{\prime},U^{\prime})$, so $C_{0}\left(  [0,1]\right)
\otimes\mathcal{A}_{\theta}\subseteq C^{\ast}(A^{\prime},U^{\prime})$, hence
$f_{2}\in C^{\ast}(A^{\prime},U^{\prime})$. Hence our lemma is proven.
\end{proof}

\begin{theorem}
\label{elliptic}Let $\theta\in\lbrack0,1)$. Let $\varphi:z\in\mathbb{D}\mapsto
e^{2i\pi\theta}z$ where $\mathbb{D}$ is the closed unit disk in $\mathbb{C}$.
Then $C(\mathbb{D})\rtimes_{\varphi}\mathbb{Z}$ is *-isomorphic to:%
\[
\left\{  f\in C\left(  [0,1],\mathcal{A}_{\theta}\right)  :f(0)\in C^{\ast
}(U)\right\}
\]
where $\mathcal{A}_{\theta}=C^{\ast}(U,V)$ with $U,V$ unitaries such that
$VU=e^{2i\pi\theta}UV$.
\end{theorem}

\begin{proof}
Since $\mathbb{D}\backslash\left\{  0\right\}  $ is fixed by $\varphi$ and
$C(\{0\})\rtimes_{\varphi}\mathbb{Z}=C(\mathbb{T})$, we have the following
easy exact sequence:%
\[
0\longrightarrow C\left(  \mathbb{D}\backslash\left\{  0\right\}  \right)
\rtimes_{\varphi}\mathbb{Z}\longrightarrow C(\mathbb{D})\rtimes_{\varphi
}\mathbb{Z}\longrightarrow C(\{0\})\rtimes_{\varphi}\mathbb{Z}\longrightarrow
0
\]
Let $\mathfrak{B}_{\theta}=\left\{  f\in C\left(  [0,1],\mathcal{A}_{\theta
}\right)  :f(0)\in C^{\ast}(U)\right\}  $. For $f\in\mathfrak{B}_{\theta}$ and
$t\in\lbrack0,1]$ we set $\zeta_{t}(f)=f(t)$. We set $A^{\prime}:t\in
\lbrack0,1]\mapsto tV$ and $U^{\prime}:t\in\lbrack0,1]\mapsto U$ in
$\mathfrak{B}_{\theta}$. Then we immediately see that $U^{\prime\ast}%
A^{\prime}U^{\prime}=\varphi(A^{\prime})$. Hence there exists a unique
*-morphism $\psi$ from $C(\mathbb{D})\rtimes_{\varphi}\mathbb{Z}$ into
$\mathfrak{B}_{\theta}$ such that $\psi(U_{\varphi})=U^{\prime}$ and
$\psi(A_{\varphi})=A^{\prime}$. By Lemma (\ref{fieldsofAtheta}),
$\mathfrak{B}_{\theta}=C^{\ast}(A^{\prime},U^{\prime})$ so $\psi$ is a
*-epimorphism onto $\mathfrak{B}_{\theta}$.

Now, assume $\psi(a)=0$ for some $a\in C\left(  \mathbb{D}\backslash\left\{
0\right\}  \right)  \rtimes_{\varphi}\mathbb{Z}$. Let $\pi$ be a nonzero
irreducible *-representation of $C\left(  \mathbb{D}\right)  \rtimes_{\varphi
}\mathbb{Z}$ on some Hilbert space $\mathcal{H}$. First, assume $\pi=\rho
\circ\xi_{r}$ for some $\rho\in\widehat{\mathcal{A}_{\theta}}$ and some
$r\in]0,1]$. We define $\pi^{\prime}=\rho\circ\zeta_{t}$ on $\mathfrak{B}%
_{\theta}$ and check that by construction $\pi=\pi^{\prime}\circ\psi$ on the
generators, hence on $C\left(  \mathbb{D}\right)  \rtimes_{\varphi}\mathbb{Z}$
(as $\pi,\pi^{\prime}\circ\psi$ *-morphisms). Therefore, $\pi(a)=0$. Now,
assume $\pi=\nu_{\omega}\circ\xi_{0}$ for some $\omega\in\mathbb{T}$.
Similarly, by setting $\pi^{\prime}=\nu_{\omega}\circ\zeta_{0}$ we see that
$\pi(a)=\pi^{\prime}\circ\psi(a)=0$. By Proposition (\ref{CovRepEll}) we
deduce that for all (irreducible) representation $\pi$ of $C\left(
\mathbb{D}\right)  \rtimes_{\varphi}\mathbb{Z}$ we have $\pi(a)=0$, so by
definition $\left\Vert a\right\Vert _{C\left(  \mathbb{D}\right)
\rtimes_{\varphi}\mathbb{Z}}=0$ and so $a=0$. Therefore, $\psi$ is a
*-isomorphism from $C\left(  \mathbb{D}\right)  \rtimes_{\varphi}\mathbb{Z}$
onto $\mathfrak{B}_{\theta}$.
\end{proof}

In particular, if $P\in C(\mathbb{D})\rtimes_{\varphi}\mathbb{Z}$ is a
projection, then it is obviously homotopic to either $0$ or $1$ (for all
$t,s\in\lbrack0,1]$ set $Q(s,t)=P(st)$), so we deduce immediately that
$K_{0}(C(\mathbb{D})\rtimes_{\varphi}\mathbb{Z})=\mathbb{Z}$. It is also
immediate that $K_{1}\left(  C(\mathbb{D})\rtimes_{\varphi}\mathbb{Z}\right)
=\mathbb{Z}$ is generated by $t\in\lbrack0,1]\mapsto U$, as any unitary $W\in
C(\mathbb{D})\rtimes_{\varphi}\mathbb{Z}$ is homotopic to a unitary of
$C(\mathbb{T})$ by setting $W^{\prime}:s,t\in\lbrack0,1]=W(st)$ and
$K_{1}\left(  C(\mathbb{T})\right)  $ is generated by $U$. Of course, these
results confirm the more general facts of the previous section.

\bigskip As in the case of hyperbolic automorphisms, some description of the
topology of $\widehat{C(\mathbb{D})\rtimes_{\varphi}\mathbb{Z}}$ is possible.
To ease notations, we identify $\widehat{C(\mathbb{D})\rtimes_{\varphi
}\mathbb{Z}}$ and $\mathfrak{B}_{\theta}$, and for all $t\in\lbrack0,1]$ we
have $\zeta_{t}=\xi_{t}$. We first define again a partition of $\widehat
{C(\mathbb{D})\rtimes_{\varphi}\mathbb{Z}}$:

\begin{definition}
Using the notations of Proposition (\ref{CovRepEll}), we define:

\begin{itemize}
\item $\operatorname*{IrrRep}_{1}=\left\{  \pi_{\omega}=\nu_{\omega}\circ
\xi_{0}:\omega\in\mathbb{T}\right\}  $ and

\item $\operatorname*{IrrRep}_{c}=\left\{  \rho\circ\xi_{r}:\rho\in
\widehat{\mathcal{A}_{\theta}},r\in]0,1]\right\}  $,
\end{itemize}

so that we have $\operatorname*{IrrRep}_{1}\cup\operatorname*{IrrRep}%
_{c}=\widehat{C(\mathbb{D})\rtimes_{\varphi}\mathbb{Z}}$ by Proposition
(\ref{CovRepEll}).
\end{definition}

We distinguish two cases: when $\theta$ is irrational, $\widehat
{C(\mathbb{D})\rtimes_{\varphi}\mathbb{Z}}$ is not Hausdorff, yet when
$\theta\in\mathbb{Q}$ then $\widehat{C(\mathbb{D})\rtimes_{\varphi}\mathbb{Z}%
}$ is a compact Hausdorff space.

\begin{proposition}
\label{ellipticspectrumirr}Let $\theta$ be irrational and let $\varphi$ be the
rotation $z\mapsto\exp\left(  2i\pi\theta\right)  z$. Define the map $q$ by
$\rho\circ\xi_{r}\in\operatorname*{IrrRep}_{c}\mapsto\left(  \rho,r\right)
\in\widehat{\mathcal{A}_{\theta}}\times]0,1]$ and $\pi_{\theta}\in
\operatorname*{IrrRep}_{1}\mapsto\omega\in\mathbb{T}$ and identify
$\operatorname*{IrrRep}$ with the spectrum of $C(\mathbb{D})\rtimes_{\varphi
}\mathbb{Z}$ as a set.\\ Then $q$ is an homeomorphism from $\widehat
{C(\mathbb{D})\rtimes_{\varphi}\mathbb{Z}}$ onto $\mathbb{T}\cup\left(
\widehat{\mathcal{A}_{\theta}}\times]0,1]\right)  $ when $\mathbb{T}%
\cup\left(  \widehat{\mathcal{A}_{\theta}}\times]0,1]\right)  $ is endowed
with the topology where closed sets are either of the form $\mathbb{T\cup
}\left(  \widehat{\mathcal{A}_{\theta}}\times F\right)  $ or $Y\cup\left(
\widehat{\mathcal{A}_{\theta}}\times X\right)  $ where $F$ is relatively
closed in $]0,1]$,the set $X$ is compact in $]0,1]$ and $Y$ is closed in
$\mathbb{T}$.
\end{proposition}

\begin{proof}
Let $F$ be a subset of $\widehat{C(\mathbb{D})\rtimes_{\varphi}\mathbb{Z}}$.
We set
\begin{equation}
T=\left\{  r\in]0,1]:\exists\rho\in\widehat{\mathcal{A}_{\theta}}%
\ \ \ \rho\circ\xi_{r}\in F\right\}  \text{.} \label{T}%
\end{equation}
Let $y$ be in the closure of $T$ in $[0,1]$. Now if $f\in C(\mathbb{D}%
)\rtimes_{\varphi}\mathbb{Z}$ such that $\pi(f)=0$ for all $\pi\in
F\cap\operatorname*{IrrRep}_{c}$, then $f(T)=\left\{  0\right\}  $, as all the
representations of $\mathcal{A}_{\theta}$ are faithful (as $\theta$ irrational
so $\mathcal{A}_{\theta}$ is simple), therefore $f(y)=0$ by continuity. In
particular, if $0$ is in the closure of $T$ then $\pi_{\omega}(f)=0$ for all
$\omega\in\mathbb{T}$ and thus $\operatorname*{IrrRep}_{1}\subseteq F$.
Moreover, since $\rho\circ\xi_{y}(f)=0$ for all $f\in%
{\displaystyle\bigcap\limits_{\pi\in F}}
\ker\pi$ and for all $\rho\in\widehat{\mathcal{A}_{\theta}}$, if $F$ is closed
then $\rho\circ\xi_{y}\in F$ for all $\rho\in\widehat{\mathcal{A}_{\theta}}$
and thus $y\in T$, so $T$ is closed and $F\cap\operatorname*{IrrRep}_{c}%
=T\cap]0,1]\times\widehat{\mathcal{A}_{\theta}}$. By \cite[Proposition 3.2.1,
p. 61]{Dixmier}, the set $\operatorname*{IrrRep}_{1}$ is closed and
homeomorphic to $\mathbb{T}$ via the map $\omega\in\mathbb{T}\mapsto
\pi_{\omega}\in\operatorname*{IrrRep}_{1}$, so if $F$ is closed then
$F\cap\operatorname*{IrrRep}_{1}$ is closed as well, and so is $Y$ in
$\mathbb{T}$. This proves the necessary conditions of our Proposition (note
that $T$ is compact in $]0,1]$ if and only if it is closed in $]0,1]$ and $0$
is not in its closure in $[0,1]$).

We now turn to the sufficient condition. Let again $F$ be a given set in
$\widehat{C(\mathbb{D})\rtimes_{\varphi}\mathbb{Z}}$ and assume that $T$, as
defined in equation (\ref{T}), is closed in $]0,1]$ and $Y$ closed in
$\mathbb{T}$. Now, set:
\[
\mathcal{M}=\left\{  f\in C(\mathbb{D})\rtimes_{\varphi}\mathbb{Z}%
:f(T)=\left\{  0\right\}  \text{ and }f(0)(Y)=\left\{  0\right\}  \right\}
\text{.}%
\]
Let $f\in\mathcal{M}$. By construction, $\pi(\mathcal{M})=\left\{  0\right\}
$ for all $\pi\in F$. Conversely, let $\pi\in\widehat{C(\mathbb{D}%
)\rtimes_{\varphi}\mathbb{Z}}$ such that $\pi(\mathcal{M})=\left\{  0\right\}
$. Now, if $\pi\in\operatorname*{IrrRep}_{1}$, then let $\mu\in\mathbb{T}$
such that $\pi=\pi_{\mu}$. Then $\pi_{\mu}(f)=f(0)(\mu)=0$ for \ all $f(0)\in
C(\mathbb{T})$ such that $f(Y)=\left\{  0\right\}  $. Thus, $\mu\in
\overline{Y}=Y$. Hence $\pi\in F$. Now, assume $\pi=\rho\circ\xi_{x}$ for
$x\in]0,1]$. Then $\pi(f)=0$ by construction. If $\tau$ is the canonical
tracial trace of $\mathcal{A}_{\theta}$, then it is easily seen that
$\widehat{\tau}(f):t\mapsto\tau(f(t))$ is a continuous surjection from
$C_{0}\left(  ]0,1],\mathcal{A}_{\theta}\right)  $ onto $C_{0}\left(
]0,1]\right)  $. Now, $f\in\mathcal{M}$ if and only if $\widehat{\tau
}(f)(t)=0$ for all $t\in T$. Therefore, if $\widehat{\tau}(f)(x)=0$ for all
$f\in\mathcal{M}$ then $x\in\overline{T}\cap]0,1]$. By assumption, $\left\{
\rho\in\widehat{\mathcal{A}_{\theta}}:\rho\circ\xi_{x}\in F\right\}
=\widehat{\mathcal{A}_{\theta}}$ (we assumed it was nonempty), so it follows
that $\pi\in F$. Hence $F=\left\{  \pi\in\widehat{C(\mathbb{D})\rtimes
_{\varphi}\mathbb{Z}}:\pi(\mathcal{M})=\left\{  0\right\}  \right\}  $, so $F$
is closed by \cite[Proposition 3.1.2, p. 60]{Dixmier}. This concludes our proof.
\end{proof}

\begin{remark}
The set $\mathbb{T}\times]0,1]$ is not in bijection with $\widehat
{C(\mathbb{D})\rtimes_{\varphi}\mathbb{Z}}$, as $\widehat{\mathcal{A}_{\theta
}}$ has many elements, though its topology is the coarse topology.
\end{remark}

\begin{proposition}
Let $\theta\in\mathbb{Q}$. Then the spectrum of $C(\mathbb{D})\rtimes
_{\varphi}\mathbb{Z}$ is homeomorphic to $\mathfrak{S}=\left\{  \left(
t\alpha,\beta\right)  :\alpha,\beta\in\mathbb{T},t\in\lbrack0,1]\right\}  $
with its usual topology.
\end{proposition}

\begin{proof}
Let $\theta=\frac{p}{q}$ with $p,q$ relatively prime, and $\omega=\exp\left(
2i\pi\theta\right)  $. Let $\Xi_{q}=C(\mathfrak{S},M_{q}(\mathbb{C}))$. For
all $\left(  \eta,\lambda\right)  \in\mathfrak{S}$ we set:
\[
U^{\prime}(\eta,\lambda)=\eta\left[
\begin{array}
[c]{cccc}%
0 & \cdots & 0 & 1\\
1 & \ddots &  & 0\\
0 & \ddots & \ddots & \vdots\\
0 & 0 & 1 & 0
\end{array}
\right]  \text{, }A^{\prime}(\eta,\lambda)=\left[
\begin{array}
[c]{cccc}%
\omega &  &  & \\
& \omega\lambda &  & \\
&  & \ddots & \\
&  &  & \omega^{q-1}\lambda
\end{array}
\right]  \text{.}%
\]
Of course, $U^{\prime}$,$A^{\prime}\in\Xi_{q}$ and $U^{\prime\ast}A^{\prime
}U^{\prime}=\omega A^{\prime}=\varphi(A^{\prime})$, and $\sigma(A^{\prime
})=\mathbb{D}$. Therefore, there exists a *-epimorphism $\rho$ from
$C(\mathbb{D})\rtimes_{\varphi}\mathbb{Z}$ onto $C^{\ast}(A^{\prime}%
,U^{\prime})$ such that $\rho(A_{\varphi})=A^{\prime}$ and $\rho(U_{\varphi
})=U^{\prime}$. Using the method of the proof of Proposition (\ref{elliptic}),
we see that $\rho$ is also injective. Hence $\widehat{C(\mathbb{D}%
)\rtimes_{\varphi}\mathbb{Z}}$ is homeomorphic to the spectrum of $C^{\ast
}(A^{\prime},U^{\prime})$. Note that $C^{\ast}(A^{\prime},U^{\prime})$ is not
$\Xi_{q}$ since $f\in C^{\ast}(A^{\prime},U^{\prime})$ implies that $f\left(
\eta,\varphi(\lambda)\right)  =\omega f(\eta,\lambda)$ for all $\left(
\eta,\lambda\right)  \in\mathfrak{S}$. Yet, it is easy to prove that
$\widehat{C^{\ast}(A^{\prime},U^{\prime})}=\mathfrak{S}$. Given $\rho
\in\widehat{\mathcal{A}_{\theta}}$ there exists $\left(  \eta,\beta\right)
\in\mathbb{T}^{2}$ such that $\rho=\pi_{\eta,\beta}$ where $\pi_{\eta,\beta
}(U)=U^{\prime}(\eta,\beta)$ and $\pi_{\eta,\beta}(V)=A^{\prime}(\eta,\beta)$,
so $\pi_{\eta,\beta}(A_{\varphi}(t))=A^{\prime}(\eta,t\beta)$. Let us define
$\rho_{\eta,t\beta}\in\widehat{C^{\ast}(A^{\prime},U^{\prime})}$ by
$\rho_{\eta,t\beta}=\pi_{\eta,\beta}\circ\zeta_{t}$ for all $t\in\lbrack0,1]$
and $\left(  \eta,\beta\right)  \in\mathbb{T}^{2}$. We then see by
\ref{CovRepEll} that $\widehat{C^{\ast}(A^{\prime},U^{\prime})}=\left\{
\rho_{x}:x\in\mathfrak{S}\right\}  $ and the map $q:\rho_{x}\mapsto
x\in\mathfrak{S}$ is an obvious bijection. Now, let $F$ be a subset of
$\widehat{C^{\ast}(A^{\prime},U^{\prime})}$ and let $M=%
{\displaystyle\bigcap\limits_{\pi\in F}}
\ker\pi$. Let $G=q(F)$ and let $f\in M$. Then $y\in\overline{G}$ if and only
if $f(y)=0$ for all $f\in M$, and thus $M\subseteq\ker\rho_{y}$. Thus, if $F$
is closed then $\rho_{y}\in F$ so $y\in G$ and thus $G$ is closed; if $G$ is
closed then $y\in G$ so $\rho_{y}\in F$ and $F=\left\{  \rho:\rho(M)=\left\{
0\right\}  \right\}  $ so $F$ is closed by \cite[Proposition 3.1.2,
p.60]{Dixmier}. Hence $q$ is an homeomorphism.
\end{proof}

\section{Appendix: Some calculations}

We group in this appendix a few details about proofs of theorems in the main
text. We include them in this version of our paper for the reader's convenience.

\begin{proposition}
An automorphism $\varphi:z\in\mathbb{D}\mapsto e^{2i\pi\theta}\frac{z-z_{0}%
}{z-\overline{z_{0}}z}$ ($\left\vert z_{0}\right\vert <1$, $\theta\in
\lbrack0,1)$), $\varphi\not =\operatorname*{Id}$, is hyperbolic, elliptic or
parabolic if and only if, respectively, $\left\vert \sin\left(  \pi
\theta\right)  \right\vert <\left\vert z_{0}\right\vert $, $\left\vert
\sin\left(  \pi\theta\right)  \right\vert >\left\vert z_{0}\right\vert $ or
$\left\vert \sin\left(  \pi\theta\right)  \right\vert =\left\vert
z_{0}\right\vert $.
\end{proposition}

\begin{proof}
We recall from Theorem (\ref{autodisk2}) that the fixed point of $\varphi$ are
the roots of $\overline{z_{0}}z^{2}+(e^{2i\pi\theta}-1)z-e^{2i\pi\theta}z_{0}%
$. This polynomial is of degree $1$ if and only if $z_{0}=0$, and in this case
$\varphi$ is a rotation (hence elliptic, unless $\theta=0$ and then
$\varphi=\operatorname*{Id}$), and we have indeed $\left\vert z_{0}\right\vert
<\left\vert \sin\left(  \pi\theta\right)  \right\vert $ for $\varphi
\not =\operatorname*{Id}$ and $z_{0}=0$. Otherwise, the discriminant of this
trinomial is $\Delta=\left(  e^{2i\pi\theta}-1\right)  ^{2}+4\lambda\left\vert
z_{0}\right\vert ^{2}$, while the modulus of the product of the roots is
$\left\vert e^{2i\pi\theta}\frac{z_{0}}{\overline{z_{0}}}\right\vert =1$,
hence either both fixed points (counting multiplicity) are on $\mathbb{T}$, or
only one is in $\mathbb{D}$, and then it lies in the interior of $\mathbb{D}$.
Since $e^{2i\pi\theta}-1=2i\exp\left(  i\pi\theta\right)  \left(  \sin\left(
\pi\theta\right)  \right)  $, we have $\Delta=4\exp\left(  2i\pi\theta\right)
\left(  \left\vert z_{0}\right\vert ^{2}-\sin^{2}\left(  \pi\theta\right)
\right)  $. Hence, $\Delta=0$ if and only if $\left\vert z_{0}\right\vert
=\left\vert \sin\left(  \pi\theta\right)  \right\vert $. Now, we denote
$\sqrt[2]{\Delta}$ for any square root of the complex $\Delta$ and $z_{\pm}$
the fixed points of $\varphi$, so that
\[
z_{\pm}\overline{z_{0}}=\frac{\left(  1-e^{2i\pi\theta}\right)  \pm
\sqrt[2]{\Delta}}{2}=e^{i\pi\theta}\left(  i\sin\left(  \pi\theta\right)
\pm\sqrt[2]{\left\vert z_{0}\right\vert ^{2}-\sin^{2}\left(  \pi\theta\right)
}\right)  \text{.}%
\]
If $\left\vert z_{0}\right\vert >\left\vert \sin\left(  \pi\theta\right)
\right\vert $ then the roots $\pm\sqrt[2]{\Delta}$of $\Delta$ are real and
distinct, so
\[
\left\vert z_{\pm}\right\vert ^{2}=\left\vert z_{0}\right\vert ^{-2}\left\vert
z_{\pm}z_{0}\right\vert ^{2}=\left\vert z_{0}\right\vert ^{-2}\left(  \left(
\sin\left(  \pi\theta\right)  \right)  ^{2}+\left\vert z_{0}\right\vert
^{2}-\sin\left(  \pi\theta\right)  ^{2}\right)  =1\text{.}%
\]
Hence, $z_{\pm}\in\mathbb{T}$ and are distinct. If $\left\vert z_{0}%
\right\vert <\left\vert \sin\left(  \pi\theta\right)  \right\vert $ then
$\left\vert \overline{z_{0}}^{-1}\left(  \lambda-1\right)  \right\vert
=2\left\vert z_{0}\right\vert ^{-1}\left\vert \sin\left(  \pi\theta\right)
\right\vert >2$. Yet if both fixed points of $\varphi$ are in $\mathbb{T}$
then the modulus of their sum is bounded above by 2. Hence, $\varphi$ has one
interior fixed point in $\mathbb{D}$ (and one outside of $\mathbb{D}$).
\end{proof}

\begin{proposition}
Let $\varphi$ be a hyperbolic automorphism of $\mathbb{D}$. Then there exists
a unique $a\in(0,1)\ $such that $\varphi$ is conformally conjugated to
$z\mapsto\frac{z+a}{1+az}$.
\end{proposition}

\begin{proof}
By Theorem (\ref{autodisk2}), there exists such an $a\in(0,1)$. To prove it is
unique, we consider the following. Let $\varphi:z\mapsto\frac{z+x_{0}}%
{1+x_{0}z}$ and $\psi:z\mapsto\frac{z+y_{0}}{1+y_{0}z}$ for $x_{0},y_{0}%
\in(0,1)$. Assume that there exists $\phi\in\mathbb{M}_{\mathbb{D}}$ such that
$\phi\circ\psi\circ\phi^{-1}=\varphi$ and define $c\in\mathbb{D}$ such that
$0=\phi(c)$. Note that $\left\vert c\right\vert <1$. Now, $\phi\circ\psi
\circ\phi^{-1}=\varphi$ implies that $\phi\left(  \left\{  -1,1\right\}
\right)  =\left\{  -1,1\right\}  $. If $\phi(1)=1$ then $\phi=z\mapsto
\frac{z-c}{1-cz}$ and thus $\phi$ commutes with $\varphi,\psi$ and thus
$\varphi=\psi$. If $\phi(1)=-1$, we check that $\phi=z\mapsto\lambda\frac
{z-c}{1-\overline{c}z}$ for some $\lambda\in\mathbb{T}$. Since $\phi(1)=-1$ we
conclude $\lambda=\frac{1-\overline{c}}{c-1}$, and since $\phi(-1)=1$ we have
$\lambda=-\frac{1+\overline{c}}{1+c}$, so $c=\overline{c}$ and thus
$\lambda=-1$. Note that we then have $\phi^{-1}=\phi$. Since $\phi\circ
\psi\circ\phi=\varphi$ we must have $\varphi(c)=\phi(y_{0})$, hence
$\frac{c+x_{0}}{1+x_{0}c}=-\frac{y_{0}-c}{1-cy_{0}}$, so $0=x_{0}c^{2}%
+cx_{0}y_{0}-y_{0}-x_{0}$. This is equivalent to $c\in\left\{  -1,1\right\}  $
or $y_{0}=-x_{0}$, yet both are impossible by assumption. Hence $x_{0}=y_{0}$.
\end{proof}

\begin{proposition}
Any two hyperbolic automorphisms of $\mathbb{D}$ are topologically conjugated.
\end{proposition}

\begin{proof}
We let $\varphi$ and $\psi$ be two hyperbolic automorphisms with attractive
fixed point $1$ and repulsive fixed point $-1$, which is sufficient by
Theorem\ (\ref{autodisk2}).

We define for all $\infty\geq r\geq1$ the circle $\Gamma_{r}$ of radius $r$
and passing by the points $-1$ and $1$ and such that $\Gamma_{r}\cap
\mathbb{D}$ is contained in the upper half plane. We denote by $\Gamma_{-r}$
the symmetric of $\Gamma_{r}$ with respect to the real axis. Note that
$\varphi(\Gamma_{\varepsilon r})=\psi(\Gamma_{\varepsilon r})=\Gamma
_{\varepsilon r}$ for all $\varepsilon\in\left\{  -1,1\right\}  $ and
$r\in\lbrack1,\infty]$. Let $\gamma_{\varphi}^{ix}$ be the arc of
$\Gamma_{r^{ix}}$ from $ix$ to $\varphi(ix)$ and let $\gamma_{\psi}^{ix}$ be
the arc of $\Gamma_{r^{ix}}$ from $ix$ to $\psi(ix),$ where all the arcs we
consider in this proof are the ones contained in $\mathbb{D}$. Let
$s_{\varphi}$ and $s_{\psi}$ be the respective length of $\gamma_{\varphi
}^{ix}$ and $\gamma_{\psi}^{ix}$ when $\Gamma_{r^{ix}}$ inherits its metric
from the usual metric on $\mathbb{C}$. Let $\omega\in\gamma_{\varphi}^{ix}$
and let $s_{\omega}$ be the length of the arc of $\Gamma_{r^{ix}}$ from $ix$
to $\omega$. We set $\mu_{0}(\omega)$ to be the unique point in $\Gamma
_{r^{ix}}$ and the right half plane so that the arc of $\Gamma_{r}$ from $ix$
to $\mu_{0}(\omega)$ has length $\frac{s_{\omega}}{s_{\varphi}}s_{\psi}$. Then
it is easy to check that $\mu_{0}(\omega)\in\gamma_{\psi}^{ix}$ and $\mu
_{0}\circ\varphi(ix)=\psi\circ\mu_{0}(ix)$ by construction. Now, let
$\mathcal{D}_{\varphi}=\cup_{x\in\lbrack-1,1]}\left(  \gamma_{\varphi}%
^{ix}\backslash\left\{  \varphi(ix)\right\}  \right)  $. By construction, for
all $z\in\mathbb{D}$, there exists a unique $y\in\mathcal{D}_{\varphi}$ such
that $\mathcal{O}_{\varphi}(z)=\mathcal{O}_{\varphi}(y)$. With obvious
notations, we therefore have constructed a bijection $\mu_{0}:\overline
{\mathcal{D}_{\varphi}}\rightarrow\overline{\mathcal{D}_{\psi}}$ which is
continuous as well (it is obvious by construction on each arc $\gamma
_{\varphi}^{ix}$. Some easy calculation shows that $\mu_{0}$ is in fact
continuous on $\overline{\mathcal{D}_{\varphi}}$). Now, since our construction
is symmetric in $\varphi$ and $\psi$, we also obtain a continuous bijection
$\mu_{0}^{\prime}:\overline{\mathcal{D}_{\psi}}\rightarrow\overline
{\mathcal{D}_{\varphi}}$ which also satisfies $\mu_{0}^{\prime}\circ\mu
_{0}=\mu_{0}^{\prime}\circ\mu_{0}=\operatorname*{Id}_{\overline{\mathcal{D}%
_{\varphi}}}$ by symmetry. Hence $\mu_{0}$ is a homeomorphism from
$\overline{\mathcal{D}_{\varphi}}$ onto $\overline{\mathcal{D}_{\psi}}$ such
that $\mu_{0}\circ\varphi=\psi\circ\mu_{0}$ on $i[-1,1]$.

By induction, we now extend $\mu_{0}$ to a homeomorphism of $\mathbb{D}$ which
conjugates $\varphi$ and $\psi$. For all $n\in\mathbb{Z}$, we set
$X_{n}^{\varphi}=\varphi^{n}\left(  \overline{\mathcal{D}_{\varphi}}\right)  $
and $X_{n}^{\psi}=\psi^{n}\left(  \overline{\mathcal{D}_{\psi}}\right)  $, as
well as $L_{n}^{\varphi}=\varphi^{n}\left(  L_{0}^{\varphi}\right)  $. Then we
assume that for some $n\in\mathbb{N}$, we have constructed the homeomorphism
$\mu_{n}$ from $X_{n}^{\varphi}$ onto $X_{n}^{\psi}$ such that $\mu_{n}%
\circ\varphi=\psi\circ\mu_{n}$ on $L_{n}^{\varphi}$. We define $\mu_{n+1}%
=\psi\circ\mu_{n}\circ\varphi^{-1}$: thus $\mu_{n+1}$ is by construction an
homeomorphism from $X_{n+1}^{\varphi}$ onto $X_{n+1}^{\psi}$. Let $x\in
L_{n+1}^{\varphi}$. We then have $\mu_{n+1}(x)=\psi\circ\mu_{n}\circ
\varphi^{-1}(x)=\mu_{n}(x)$ by our induction hypothesis. Hence $\mu_{n}$ and
$\mu_{n+1}$ agree on $L_{n+1}^{\varphi}$. Hence $\mu_{n+1}\left(
\varphi(x)\right)  =\psi\left(  \mu_{n}(x)\right)  =\psi\circ\mu_{n+1}(x)$,
hence our induction is complete since we have already constructed $\mu_{0}$.
The same construction yields $\mu_{n}$ for all $n\leq0$ as well.

We now define the map $\mu$ by setting $\mu(x)=\mu_{n}(x)$ for any
$x\in\mathbb{D}\backslash\left\{  -1,1\right\}  $ where $n\in\mathbb{Z}$ is
defined by $x\in X_{n}^{\varphi}$. Now, $\mu$ is well-defined since if $x\in
X_{n}^{\varphi}\cap X_{n+1}^{\varphi}=L_{n+1}^{\varphi}$ then $\mu
_{n+1}(x)=\mu_{n}(x)$ by the previous induction. It is now immediate that
$\mu$ is an homeomorphism from $\mathbb{D}\backslash\left\{  -1,1\right\}  $
onto itself. We also set $\mu(\varepsilon)=\varepsilon$ for $\varepsilon
\in\left\{  -1,1\right\}  $. We now check that $\mu$ is also continuous at $1
$ and $-1$. By symmetry, it suffices to prove continuity at $1$. Let $\left(
x_{n}\right)  _{n\in\mathbb{N}}$ be a sequence in $\mathbb{D}$ with limit $1$.
Let $\Omega$ be a neighborhood of $1$ in $\mathbb{D}$. Then there exists
$M\in\mathbb{N}$ such that $\cup_{j\geq M}X_{j}^{\psi}\subset\Omega$. Now
there exists $N\in\mathbb{N}$ such that $x_{n}\in\cup_{j\geq M}X_{j}^{\varphi
}$ for all $n\geq N$. Hence $\mu(x_{n})\in\cup_{j\geq M}X_{j}^{\psi}%
\subset\Omega$ for all $n\geq N$ and thus $\lim_{n\rightarrow\infty}\mu
(x_{n})=1$.

By symmetry of $\varphi$ and $\psi$ again, $\mu$ is a homeomorphism of
$\mathbb{D}$. Now, let $x\in\mathbb{D}\backslash\{-1,1\}$ and let
$n\in\mathbb{Z}$ such that $x\in X_{n}^{\varphi}$. Then $\mu\circ
\varphi(x)=\mu_{n+1}\left(  \varphi(x)\right)  =\psi\left(  \mu_{n+1}%
(x)\right)  =\psi\circ\mu(x),$ hence $\mu$ conjugates $\psi$ and $\varphi$ (by
construction $\mu\circ\varphi(\varepsilon)=\psi\circ\mu(\varepsilon)$ when
$\varepsilon\in\left\{  -1,1\right\}  $).
\end{proof}

\bibliographystyle{amsplain}
\bibliography{thesis}

\end{document}